# NUMERICAL APPROXIMATIONS FOR THE BINARY FLUID-SURFACTANT PHASE FIELD MODEL WITH FLUID FLOW: SECOND-ORDER, LINEAR, ENERGY STABLE SCHEMES

XIAOFENG YANG

ABSTRACT. In this paper, we consider numerical approximations of a binary fluid-surfactant phase-field model coupled with the fluid flow, in which the system is highly nonlinear that couples the incompressible Navier-Stokes equations and two Cahn-Hilliard type equations. We develop two, linear and second order time marching schemes for solving this system, by combining the "Invariant Energy Quadratization" approach for the nonlinear potentials, the projection method for the Navier-Stokes equation, and a subtle implicit-explicit treatment for the stress and convective terms. We prove the well-posedness of the linear system and its unconditionally energy stability rigorously. Various 2D and 3D numerical experiments are performed to validate the accuracy and energy stability of the proposed schemes.

## 1. INTRODUCTION

Surfactant is a substance that tends to reduce the surface tension of a liquid in which it is dissolved, and allows for the mixing of multiple immiscible liquids. Thus the modeling and numerical simulations for investigating the binary fluid-surfactant system naturally incorporate the hydrodynamics and the free moving interfaces. In the pioneering work of Laradji et. al. [21, 22], the phase field method, an effective modeling and simulation tool to resolve the motion of free interfaces, was firstly used to study the phase transition behaviors of the monolayer microemulsion system formed by soluble surfactant molecules. Since then, a variety of binary fluid surfactant phase field models had been developed in the past two decades, see [21, 22, 36, 20, 3, 8, 35]. Concerning the widely applications of phase field methods in science and engineering and some recent developments in the related computational technologies, we refer to [17, 24, 9, 39, 4, 27, 34, 28, 19, 25, 10] and literatures cited therein.

We now give a brief review about the commonly used soluble fluid surfactant models using the phase field approach. In [21, 22], two phase field variables are introduced to represent the local *densities* of the fluids, and the local *concentration* of the surfactant. The model is then formulated by the variational approach of the total free energy where two types of nonlinear energy terms are considered, including the phenomenological Ginzburg-Landau double well potential for the density variable to describe the phase separation behaviors of the fluid mixture; and the nonlinear coupling entropy term to ensure the high fraction of the surfactant near the fluid interface. Following the work of Laradji, the use of two phase variables for surfactant models formed the general structure for the surfactant phase field

Department of Mathematics, University of South Carolina, Columbia, SC 29208, USA. Email: xfyang@math.sc.edu. X. Yang's research is partially supported by the U.S. National Science Foundation under grant numbers DMS-1720212 and DMS-1418898.

models. In [20], the authors modified Laradji's model by adding an extra diffusion term and a Ginzburg Landau type potential for the concentration variable, in order to improve the stability. In [36], the logarithmic Flory-Huggins potential was used to restrict the range of the concentration variable, while the nonlinear coupling entropy is essentially the same as that in [21, 22, 20]. A slightly different nonlinear coupling entropy was presented in [8], which could penalize the concentration to accumulate along the fluid interface. In [35], the authors further modified the model in [8] by adding the Flory-Huggins potential for the local concentration variable as well. Based on these modeling contributions, a large quantity of theoretical/numerical works had been devoted either to analyze the PDE system [37, 1, 7], or to perform computer simulations [11, 35, 18], or to improve the model for realistic surfactant systems [10, 3, 45, 26], etc. We refer to [23] for a thorough review about these literatures.

While numerous fluid surfactant phase models had been given attentions for more than twenty years, a non-negligible fact is almost all developed numerical algorithms had been focused on the *partial* model, i.e., no flow field had been considered, which implies that the inherent numerical challenges are greatly diminished. For the partial model, one only needs to consider how to discretize the nonlinear stiff terms, for instance, the cubic term induced by the double well potential, or the term induced by the logarithmic Flory-Huggins potential. To name a few, in [18], a second order in time, nonlinear scheme was developed based on the Crank-Nicolson method to solve the partial model, however, the scheme is not provably energy stable. In [11], the authors developed a first order in time, nonlinear scheme based on the convex splitting approach, where the convex part of the free energy potential is treated implicitly while the concave part is treated explicitly, but the arguments about the convex-concave decomposition for the coupling potential are not valid as well since it is not sufficient to justify the convexity of a function with multiple variables from the positivity of second order partial derivatives. In addition, their scheme is only of first order in time, and its computational cost is relatively expensive due to the nonlinear nature. In [42], the authors developed two provably, unconditionally energy stable, linear and second order schemes. But the model considered in [42] is still the partial model without considering the flow filed.

Compared with the partial model, in addition to the stiffness issue mentioned above, the full hydrodynamical surfactant model is conceivably more complicated for algorithms design since one must consider to develop efficient temporal discretizations for the nonlinear coupling terms between the flow field and phase field variables through stress and convective terms. Simple fully-implicit or explicit type discretizations for those terms will induce unfavorable nonlinear and/or energy unstable schemes, so they are not efficient in practice. Therefore, in this paper, the main purpose is to develop unconditionally energy stable schemes for solving the full hydrodynamics coupled fluid-surfactant model. Specifically, we choose the model developed in [8, 35] since it is not only a typical representative of nonlinearly coupled multivariate fluid-surfactant models, but also *mathematically reasonable*, i.e., the total energy is naturally bounded from below.

We emphasize that our aim of algorithm designs, i.e., the "unconditional energy stable schemes", only means the schemes have no constraints on the time step size from stability point of view. However, large time step size will definitely induce large errors in practice. This fact motivates us to develop more accurate schemes, e.g., the second order time stepping schemes while preserving the unconditional energy stability in this paper. Overall, we expect that our schemes can combine the following three desired properties, i.e., (i) *accurate* (second order in time); (ii) *stable* (the unconditional energy dissipation law holds); and (iii) *easy to implement and efficient* (only need to solve some fully *linear* equations at each time step). To achieve such a goal, by combining some well-known approaches like the projection method for Navier-Stokes equations, the Invariant Energy Quadratizaiton (IEQ) approach for nonlinear stiff terms, and implicit-explicit treatments for nonlinear stresses, we obtain two linear, second order, and provably unconditionally energy stable time marching schemes. Moreover, we rigorously prove that the well-posedness and unconditional energy stabilities hold, and demonstrate the stability and the accuracy of the proposed schemes through a number of classical benchmark phase separation simulations. To the best of the author's knowledge, the proposed schemes here are the first linear schemes with second order temporal accuracy for the hydrodynamics coupled phase field surfactant model with provable unconditional energy stabilities.

The rest of the paper is organized as follows. In Section 2, we introduce the binary fluid-surfactant phase field model developed in [8, 35] and reformulate it to couple the Navier-Stokes equations. In Section 3, we present two second order numerical schemes for simulating the target model, and rigorously prove that the produced linear systems are well-posed and the schemes satisfy the unconditional energy stabilities. Various 2D and 3D numerical experiments are carried out in Section 4 to verify the accuracy and stability of these schemes. Finally, some concluding remarks are given in Section 5.

## 2. The flow coupled surfactant model and its energy law

We fix some notations first. We consider an open, bounded, connected domain $\Omega \in \mathbb{R}^d, d = 2, 3$ with $\partial\Omega$ Lipschitz continuous. For any $f(\boldsymbol{x}), g(\boldsymbol{x}) \in L^2(\Omega)$, we denote the inner product and $L^2$ norm as

$$(f, g) = \int_\Omega f(\boldsymbol{x})g(\boldsymbol{x})d\boldsymbol{x}, \quad \|f\| = \int_\Omega |f(\boldsymbol{x})|^2 d\boldsymbol{x}. \tag{2.1}$$

Let $W^{k,r}(\Omega)$ stands for the standard Sobolev spaces equipped with the standard Sobolev norms $\|\cdot\|_{k,r}$. For $r = 1$, we write $H^1(\Omega)$ for $W^{1,2}(\Omega)$ and its corresponding norm is $\|\cdot\|_{H^1}$.

We now give a brief description for the total free energy for the fluid-surfactant phase field model developed in [8, 35]. In the binary fluid-surfactant system, monolayers of surfactant molecules form microemulsions as a random phase. Such microemulsion system usually exhibits various interesting microstructures, depending on the temperature or the composition. The phase field modeling approach simulates the dynamics of microphase separation in microemulsion systems using two phase field variables (order parameters).

To label the local densities of the two fluids, fluid I and fluid II (such as water and oil), a phase field variable $\phi(\boldsymbol{x}, t)$ is introduced such that

$$\phi(\boldsymbol{x}, t) = \begin{cases} -1 & \text{fluid I}, \\ 1 & \text{fluid II}, \end{cases} \tag{2.2}$$

with a thin smooth transition layer of width $O(\epsilon)$ connecting the two fluids. Thus the interface of the mixture is described by the zero level set $\Gamma_t = \{\boldsymbol{x} : \phi(\boldsymbol{x}, t) = 0\}$. Let $F(\phi) = (\phi^2 - 1)^2$ be the classical double well potential, one can define the mixing free energy functional as

$$\boldsymbol{E}_1(\phi) = \int_\Omega \Big(\frac{\epsilon}{2}|\nabla\phi|^2 + \frac{1}{4\epsilon}F(\phi)\Big)d\boldsymbol{x}. \tag{2.3}$$

In fact, $\boldsymbol{E}_1(\phi)$ is the most commonly used free energy in phase field models so far, where the first term in (2.3) contributes to the hydrophilic type (tendency of mixing) of interactions between the materials and the second term represents the hydrophobic type (tendency of separation) of interactions. As the consequence of the competition between the two types of interactions, the equilibrium configuration will include a diffusive interface. The phase field models that are derived by this part of free energy had been extensively studied, see [14, 2, 16, 25] and the references therein.

To represent the local concentration of surfactants, another phase field variable $\rho(\boldsymbol{x}, t)$ is introduced and its associated free energy is

$$\boldsymbol{E}_2(\rho) = \int_\Omega (\frac{\eta}{2}|\nabla\rho|^2 + \beta G(\rho))d\boldsymbol{x}, \tag{2.4}$$

where $G(\rho) = \rho\ln\rho + (1-\rho)\ln(1-\rho)$, $\beta$ and $\eta$ (both are in the same scale as $\epsilon$) are two positive parameters. Here we add the gradient potential in order to enhance the stability. We note that $G(\rho)$ is the logarithmic Flory-Huggins type energy potential which restricts the value of $\rho$ to be inside the domain of $(0, 1)$, and $\rho$ will reach its upper bound if the interface is fully saturated with surfactant (cf. [36]).

Lastly, due to the special property of surfactants that can alter the interfacial tension, the fraction of the surfactant staying at the fluid interface is high. Thus a local, nonlinear coupling entropy term between $\phi$ and $\rho$ is imposed as

$$\boldsymbol{E}_3(\phi, \rho) = \frac{\alpha}{2}\int_\Omega (\rho - |\nabla\phi|)^2 d\boldsymbol{x}, \tag{2.5}$$

where $\alpha$ is a positive parameter. This part of energy is the penalty term that enables the concentration to accumulate near the interface with a relatively high value.

When coupled with the hydrodynamics, the total energy of the binary mixture fluid-surfactant system is the sum of the kinetic energy $\boldsymbol{E}_k$, the mixing energy $\boldsymbol{E}_1$, the logarithmic energy $\boldsymbol{E}_2$, and the coupling entropy term $\boldsymbol{E}_3$, as follows,

$$\boldsymbol{E}_{tot}(\boldsymbol{u}, \phi, \rho) = \boldsymbol{E}_k(\varrho, \boldsymbol{u}) + \boldsymbol{E}_1(\phi) + \boldsymbol{E}_2(\rho) + \boldsymbol{E}_3(\phi, \rho), \tag{2.6}$$

where $\boldsymbol{E}_k(\varrho, \boldsymbol{u}) = \int_\Omega \frac{1}{2}\varrho|\boldsymbol{u}|^2 d\boldsymbol{x}$, $\varrho$ is the volume-averaged density of the mixture system, $\boldsymbol{u}$ is the is the volume-averaged velocity field.

In this paper, we consider the case where the two fluid components have matching density and viscosity, i.e., $\varrho_1 = \varrho_2 = 1$ and $\nu_1 = \nu_2 = \nu$ for simplicity. Following the generalized Onsager's principle [30], we derive the following governing system of equations and present them in dimensionless forms:

$$\boldsymbol{u}_t + (\boldsymbol{u}\cdot\nabla)\boldsymbol{u} + \nabla p - \nu\Delta\boldsymbol{u} + \phi\nabla\mu_\phi + \rho\nabla\mu_\rho = 0, \tag{2.7}$$

$$\nabla\cdot\boldsymbol{u} = 0, \tag{2.8}$$

$$\phi_t + \nabla\cdot(\boldsymbol{u}\phi) = M_1\Delta\mu_\phi, \tag{2.9}$$

$$\mu_\phi = -\epsilon\Delta\phi + \frac{1}{\epsilon}\phi(\phi^2-1) + \alpha\nabla\cdot\Big((\rho - |\nabla\phi|)\boldsymbol{Z}\Big), \tag{2.10}$$

$$\rho_t + \nabla\cdot(\boldsymbol{u}\rho) = M_2\Delta\mu_\rho, \tag{2.11}$$

$$\mu_\rho = -\eta\Delta\rho + \alpha(\rho - |\nabla\phi|) + \beta\ln\Big(\frac{\rho}{1-\rho}\Big), \tag{2.12}$$

where $\boldsymbol{Z}(\phi) = \frac{\nabla\phi}{|\nabla\phi|}$, $M_1$ and $M_2$ are two mobility parameters. For simplicity, we take the periodic boundary conditions to remove all complexities from the boundary integrals.

It is then straightforward to obtain the PDE energy law for the flow coupled system (2.7)-(2.12). By taking the $L^2$ inner products of (2.7) with $\boldsymbol{u}$, of (2.9) with $\mu_\phi$, of (2.10) with $-\phi_t$, of (2.11) with $\mu_\rho$, and of (2.12) with $-\rho_t$, using integration by parts, and combining all equalities, we obtain

$$\frac{d}{dt}\boldsymbol{E}_{tot}(\boldsymbol{u},\phi,\rho) = -\nu\|\nabla\boldsymbol{u}\|^2 - M_1\|\nabla\mu_\phi\|^2 - M_2\|\nabla\mu_\rho\|^2 \le 0. \tag{2.13}$$

In the sequel, our goal is to design temporal approximation schemes which satisfy the discrete version of the continuous energy law (2.13).

## 3. Numerical schemes

We next construct schemes for discretizing the flow coupled surfactant model (2.7)-(2.12) in time. To this end, there is a number of numerical challenges, including

- (i) how to discretize the cubic term associated with the double well potential;
- (ii) how to discretize the the logarithmic term induced by the Flory-Huggins potential;
- (iii) how to discretize the local coupling entropy terms associated with $\rho$ and $\phi$;
- (iv) how to solve the Navier-Stokes equations; and
- (v) how to discretize the nonlinear convective and stress terms.

For (i), (ii) and (iii), we adopt the recently developed IEQ approach (cf. [42, 40, 47, 43, 38, 41, 15, 46, 44]) to design desired numerical schemes. The intrinsic idea of the IEQ method is to transform the nonlinear potential into quadratic forms in terms of some new variables. This method is workable since we notice that all nonlinear potentials in the free energy are bounded from below. For (iv), we use the second order pressure correction method [13], where the pressure is then decoupled from the computations of the velocity. For (v), we use a subtle implicit-explicit treatment to obtain the fully linear schemes.

Following the work in [5, 6], we regularize the Flory-Huggins potential from domain $(0,1)$ to $(-\infty,\infty)$, where the logarithmic functional is replaced by a $C^2$ continuous, convex and piecewisely defined function. More precisely, for any $\widehat{\epsilon} > 0$, the regularized Flory-Huggins potential is given by

$$\widehat{G}(\rho) = \begin{cases} \rho\ln\rho + \frac{(1-\rho)^2}{2\widehat{\epsilon}} + (1-\rho)\ln\widehat{\epsilon} - \frac{\widehat{\epsilon}}{2}, & \text{if} \quad \rho \ge 1-\widehat{\epsilon}, \\ \rho\ln\rho + (1-\rho)\ln(1-\rho), & \text{if} \quad \widehat{\epsilon} \le \rho \le 1-\widehat{\epsilon}, \\ (1-\rho)\ln(1-\rho) + \frac{\rho^2}{2\widehat{\epsilon}} + \rho\ln\widehat{\epsilon} - \frac{\widehat{\epsilon}}{2}, & \text{if} \quad \rho \le \widehat{\epsilon}. \end{cases} \tag{3.1}$$

When $\widehat{\epsilon} \to 0$, $\widehat{G}(\rho) \to G(\rho)$. It was proved in [5, 6] that the error bound between the regularized system and the original system is controlled by $\widehat{\epsilon}$ up to a constant scaling for the Cahn-Hilliard equation. Thus we consider the numerical solution to the model formulated with the regularized functional $\widehat{G}(\rho)$, but omit the $\widehat{}$ in the notation for convenience.

It is obvious that $G(\rho)$ is bounded from below although it is not always positive in the whole domain. Thus we can rewrite the free energy functional to the following equivalent form:

$$\begin{aligned}\boldsymbol{E}_{tot}(\boldsymbol{u},\phi,\rho) = \int_\Omega \Big(\frac{1}{2}|\boldsymbol{u}|^2 + \frac{\epsilon}{2}|\nabla\phi|^2 + \frac{\eta}{2}|\nabla\rho|^2 + \frac{1}{4\epsilon}(\phi^2-1)^2 \\ + \frac{\alpha}{2}(\rho - |\nabla\phi|)^2 + \beta(\sqrt{G(\rho)+A}\,)^2\Big)d\boldsymbol{x} - \beta A|\Omega|,\end{aligned} \tag{3.2}$$

where $A$ is a positive constant to ensure $G(\rho) + A > 0$ (in all numerical examples, we use $A = 1$). Note the free energy is still the same because we simply add a zero term $\beta A - \beta A$ to the energy density functional. Then we define three auxiliary variables to be the square root of the nonlinear potentials $F(\phi)$, $(\rho - |\nabla\phi|)^2$ and $G(\rho) + A$ by

$$U = \phi^2 - 1, \tag{3.3}$$

$$V = \rho - |\nabla\phi|, \tag{3.4}$$

$$W = \sqrt{G(\rho) + A}. \tag{3.5}$$

In turn, the total free energy (3.2) can be transformed as

$$\begin{aligned}\boldsymbol{E}_{tot}(\boldsymbol{u},\phi,\rho,U,V,W) = \int_\Omega \Big(\frac{1}{2}|\boldsymbol{u}|^2 + \frac{\epsilon}{2}|\nabla\phi|^2 + \frac{\eta}{2}|\nabla\rho|^2 + \frac{1}{4\epsilon}U^2 + \frac{\alpha}{2}V^2 + \beta W^2\Big)d\boldsymbol{x} \\ - \beta A|\Omega|.\end{aligned} \tag{3.6}$$

Then we obtain a new, but equivalent partial differential system as follows:

$$\boldsymbol{u}_t + (\boldsymbol{u}\cdot\nabla)\boldsymbol{u} + \nabla p - \nu\Delta\boldsymbol{u} + \phi\nabla\mu_\phi + \rho\nabla\mu_\rho = 0, \tag{3.7}$$

$$\nabla\cdot\boldsymbol{u} = 0, \tag{3.8}$$

$$\phi_t + \nabla\cdot(\phi\boldsymbol{u}) = M_1\Delta\mu_\phi, \tag{3.9}$$

$$\mu_\phi = -\epsilon\Delta\phi + \frac{1}{\epsilon}\phi U + \alpha\nabla\cdot(V\boldsymbol{Z}), \tag{3.10}$$

$$\rho_t + \nabla\cdot(\rho\boldsymbol{u}) = M_2\Delta\mu_\rho, \tag{3.11}$$

$$\mu_\rho = -\eta\Delta\rho + \alpha V + \beta H(\rho)W, \tag{3.12}$$

$$U_t = 2\phi\phi_t, \tag{3.13}$$

$$V_t = \rho_t - \boldsymbol{Z}\cdot\nabla\phi_t, \tag{3.14}$$

$$W_t = \frac{1}{2}H(\rho)\rho_t, \tag{3.15}$$

with $H(\rho) = \frac{g(\rho)}{\sqrt{G(\rho)+A}}$ where $g(\rho) = G'(\rho)$.

The initial conditions are correspondingly

$$\begin{cases}\boldsymbol{u}|_{(t=0)} = \boldsymbol{u}_0,\ \phi|_{(t=0)} = \phi_0,\ \rho|_{(t=0)} = \rho_0, \\ U|_{(t=0)} = \phi_0^2 - 1,\ V|_{(t=0)} = \rho_0 - |\nabla\phi_0|,\ W|_{(t=0)} = \sqrt{G(\rho_0)+A}.\end{cases} \tag{3.16}$$

The transformed system still follows the energy dissipation law. By taking the sum of the $L^2$ inner products of (3.7) of $\boldsymbol{u}$, of (3.9) with $\mu_\phi$, of (3.10) with $-\phi_t$, of (3.11) with $\mu_\rho$, of (3.12) with $-\rho_t$, of (3.13) with $\frac{1}{2\epsilon}U$, of (3.14) with $\alpha V$, of (3.15) with $2\beta W$, using the integration by parts, we then obtain the energy dissipation law of the new system as (3.7)-(3.15) as

$$\frac{d}{dt}\boldsymbol{E}_{tot}(\boldsymbol{u},\phi,\rho,U,V,W) = -\nu\|\nabla\boldsymbol{u}\|^2 - M_1\|\nabla\mu_\phi\|^2 - M_2\|\nabla\mu_\rho\|^2 \le 0. \tag{3.17}$$

In the following, we focus on designing numerical schemes for time stepping of the transformed system (3.7)-(3.15), that are linear and satisfy discrete analogues of the energy law (3.17). We fix some notations here. Let $\delta t > 0$ denote the time step size and set $t_n = n\,\delta t$ for $0 \le n \le N$ with the ending time $T = N\,\delta t$. We define three Sololev spaces $H^1_{per}(\Omega) = \{\phi \in H^1(\Omega) : \phi \textit{ is periodic}\}$, $\bar{H}^1(\Omega) = \{\phi \in H^1_{per}(\Omega) : \int_\Omega \phi d\boldsymbol{x} = 0\}$, $H^1_{\boldsymbol{u}}(\Omega) = \{\boldsymbol{u} \in [H^1_{per}(\Omega)]^d\}$.

### 3.1. **Adam-Bashforth Scheme.**

We first construct a second order time stepping scheme to solve the system (3.7)-(3.15), based on the Adam-Bashforth backward differentiation formulas (BDF2).

Assuming that $(\boldsymbol{u},\phi,\rho,U,V,W)^{n-1}$ and $(\boldsymbol{u},\phi,\rho,U,V,W)^{n}$ are known:

**step 1**: we update $(\tilde{\boldsymbol{u}},\phi,\rho,U,V,W)^{n+1}$ as follows,

$$\frac{3\tilde{\boldsymbol{u}}^{n+1} - 4\boldsymbol{u}^n + \boldsymbol{u}^{n-1}}{2\delta t} + B(\boldsymbol{u}^\star,\tilde{\boldsymbol{u}}^{n+1}) - \nu\Delta\tilde{\boldsymbol{u}}^{n+1} + \nabla p^n + \phi^\star\nabla\mu_\phi^{n+1} + \rho^\star\nabla\mu_\rho^{n+1} = 0, \tag{3.18}$$

$$\frac{3\phi^{n+1} - 4\phi^n + \phi^{n-1}}{2\delta t} + \nabla\cdot(\tilde{\boldsymbol{u}}^{n+1}\phi^\star) = M_1\Delta\mu_\phi^{n+1}, \tag{3.19}$$

$$\mu_\phi^{n+1} = -\epsilon\Delta\phi^{n+1} + \frac{1}{\epsilon}\phi^\star U^{n+1} + \alpha\nabla\cdot(V^{n+1}\boldsymbol{Z}^\star), \tag{3.20}$$

$$\frac{3\rho^{n+1} - 4\rho^n + \rho^{n-1}}{2\delta t} + \nabla\cdot(\tilde{\boldsymbol{u}}^{n+1}\rho^\star) = M_2\Delta\mu_\rho^{n+1}, \tag{3.21}$$

$$\mu_\rho^{n+1} = -\eta\Delta\rho^{n+1} + \alpha V^{n+1} + \beta H^\star W^{n+1}, \tag{3.22}$$

$$3U^{n+1} - 4U^n + U^{n-1} = 2\phi^\star(3\phi^{n+1} - 4\phi^n + \phi^{n-1}), \tag{3.23}$$

$$3V^{n+1} - 4V^n + V^{n-1} = (3\rho^{n+1} - 4\rho^n + \rho^{n-1}) - \boldsymbol{Z}^\star\cdot\nabla(3\phi^{n+1} - 4\phi^n + \phi^{n-1}), \tag{3.24}$$

$$3W^{n+1} - 4W^n + W^{n-1} = \frac{1}{2}H^\star(3\rho^{n+1} - 4\rho^n + \rho^{n-1}), \tag{3.25}$$

with periodic boundary condition being imposed, where

$$\begin{aligned} &B(\boldsymbol{u},\boldsymbol{v}) = (\boldsymbol{u}\cdot\nabla)\boldsymbol{v} + \frac{1}{2}(\nabla\cdot\boldsymbol{u})\nabla\boldsymbol{v}, \\ &\boldsymbol{u}^\star = 2\boldsymbol{u}^n - \boldsymbol{u}^{n-1}, \phi^\star = 2\phi^n - \phi^{n-1}, \rho^\star = 2\rho^n - \rho^{n-1}, \\ &\boldsymbol{Z}^\star = \boldsymbol{Z}(\phi^\star), H^\star = H(\rho^\star). \end{aligned} \tag{3.26}$$

**step 2**: we update $p^{n+1}$ as follows,

$$\frac{3}{2\delta t}(\boldsymbol{u}^{n+1} - \tilde{\boldsymbol{u}}^{n+1}) + \nabla(p^{n+1} - p^n) = 0, \tag{3.27}$$

$$\nabla\cdot\boldsymbol{u}^{n+1} = 0, \tag{3.28}$$

with periodic boundary conditions.

**Remark 3.1.** *Note the scheme* (3.18)-(3.25) *is purely linear, and the computations of the variables* $(\phi, \mu_\phi, \rho, \mu_\rho, \tilde{\boldsymbol{u}})^{n+1}$ *are decoupled from the pressure* $p^{n+1}$ *by combining the second order pressure correction scheme with a subtle implicit-explicit treatment for the stress and convective terms. We are not clear how to develop second order, energy stable schemes that can decouple the computations of all variables. Indeed, it is still an open problem on how to decouple the phase field variable from the velocity field* $\boldsymbol{u}$ *while preserving the unconditional energy stability and second order temporal accuracy. So far, the best, decoupled type energy stable schemes with provable energy stabilities are just first order accurate in time (cf.* [31, 32, 33, 29]*).*

**Remark 3.2.** *The adopted projection method here was analyzed in* [12] *where it is shown (discrete time, continuous space) that the scheme is second order accurate for the velocity in* $\ell^2(0,T;L^2(\Omega))$ *but only first order accurate for the pressure in* $\ell^\infty(0,T;L^2(\Omega))$ *due to the artificial homogenous Neumann boundary condition imposed on the pressure. Since we use the periodic boundary conditions for all variables, the order of accuracy of the pressure can be second order accurate, that is shown by the numerical tests in section 4.*

Note that the new variables $U$, $V$ and $W$ will not bring up extra computational cost through the following procedure. We rewrite the equations (3.23)–(3.25) as follows,

$$\begin{cases} U^{n+1} = A_1 + 2\phi^\star\phi^{n+1}, \\ V^{n+1} = B_1 + \rho^{n+1} - \boldsymbol{Z}^\star \cdot \nabla\phi^{n+1}, \\ W^{n+1} = C_1 + \dfrac{1}{2}H^\star\rho^{n+1}, \end{cases} \tag{3.29}$$

where

$$\begin{cases} A_1 = U^{n\dagger} - 2\phi^\star\phi^{n\dagger}, \\ B_1 = V^{n\dagger} - \rho^\dagger + \boldsymbol{Z}^\star \cdot \nabla\phi^{n\dagger}, \\ C_1 = W^{n\dagger} - \dfrac{1}{2}H^\star\rho^{n\dagger}, \end{cases} \tag{3.30}$$

where $S^{n\dagger} = \frac{4S^n - S^{n-1}}{3}$ for any variable $S$.

Thus the system (3.18)-(3.25) can be rewritten as

$$\tilde{\boldsymbol{u}}^{n+1} + \frac{2\delta t}{3}B(\boldsymbol{u}^\star, \tilde{\boldsymbol{u}}^{n+1}) - \frac{2\delta t}{3}\nu\Delta\tilde{\boldsymbol{u}}^{n+1} + \frac{2\delta t}{3}\phi^\star\nabla\mu_\phi^{n+1} + \frac{2\delta t}{3}\rho^\star\nabla\mu_\rho^{n+1} = \boldsymbol{u}^{n\dagger} - \frac{2\delta t}{3}\nabla p^n, \tag{3.31}$$

$$\phi^{n+1} + \frac{2\delta t}{3}\nabla\cdot(\tilde{\boldsymbol{u}}^{n+1}\phi^\star) - \frac{2\delta t}{3}M_1\Delta\mu_\phi^{n+1} = \phi^{n\dagger}, \tag{3.32}$$

$$-\mu_\phi^{n+1} - \epsilon\Delta\phi^{n+1} + \frac{1}{\epsilon}2\phi^\star\phi^\star\phi^{n+1} + \alpha\nabla\cdot(\rho^{n+1}\boldsymbol{Z}^\star - \boldsymbol{Z}^\star\cdot\nabla\phi^{n+1}\boldsymbol{Z}^\star) = -\frac{1}{\epsilon}\phi^\star A_1 - \alpha\nabla\cdot(B_1\boldsymbol{Z}^\star), \tag{3.33}$$

$$\rho^{n+1} + \frac{2\delta t}{3}\nabla\cdot(\tilde{\boldsymbol{u}}^{n+1}\rho^\star) - \frac{2\delta t}{3}M_2\Delta\mu_\rho^{n+1} = \rho^{n\dagger}, \tag{3.34}$$

$$-\mu_\rho^{n+1} - \eta\Delta\rho^{n+1} + \alpha\rho^{n+1} - \alpha\boldsymbol{Z}^\star\cdot\nabla\phi^{n+1} + \frac{1}{2}\beta H^\star H^\star\rho^{n+1} = -\alpha B_1 - \beta H^\star C_1. \tag{3.35}$$

Therefore, the practical implementation for solving the scheme (3.18)-(3.28) is solving the above system to obtain $(\tilde{\boldsymbol{u}}, \phi, \rho)^{n+1}$, and then update $(U,V,W)^{n+1}$ through (3.29).

Now we study the wellposeness of the corresponding weak form of the semi-discretized system (3.31)-(3.35).

We define

$$\bar{\phi} = \frac{1}{|\Omega|}\int_\Omega \phi d\boldsymbol{x}, \qquad \bar{\rho} = \frac{1}{|\Omega|}\int_\Omega \rho d\boldsymbol{x},$$
$$\bar{\mu}_\phi = \frac{1}{|\Omega|}\int_\Omega \mu_\phi d\boldsymbol{x}, \quad \bar{\mu}_\rho = \frac{1}{|\Omega|}\int_\Omega \mu_\rho d\boldsymbol{x}. \tag{3.36}$$

By integrating (3.19) and (3.21), we obtain

$$\bar{\phi}^{n+1} = \bar{\phi}^n = \cdots = \bar{\phi}^0, \quad \bar{\rho}^{n+1} = \bar{\rho}^n = \cdots = \bar{\rho}^0. \tag{3.37}$$

We let

$$\phi = \phi^{n+1} - \bar{\phi}^0, \mu_\phi = \mu_\phi^{n+1} - \bar{\mu}_\phi^0,$$
$$\rho = \rho^{n+1} - \bar{\rho}^0, \mu_\rho = \mu_\rho^{n+1} - \bar{\mu}_\rho^0,$$

thus $\phi, \rho, \mu_\phi, \mu_\rho \in \bar{H}^1(\Omega)$. The weak form for (3.31)-(3.35) can be written as the following system with unknowns $\phi$, $\rho$, $\mu_\phi$, $\mu_\rho \in \bar{H}^1(\Omega)$, $\boldsymbol{u} \in H^1_{\boldsymbol{u}}(\Omega)$,

$$(\boldsymbol{u}, \boldsymbol{v}) + \frac{2\delta t}{3}(B(\boldsymbol{u}^\star, \boldsymbol{u}), \boldsymbol{v}) + \frac{2\delta t}{3}\nu(\nabla \boldsymbol{u}, \nabla \boldsymbol{v}) + \frac{2\delta t}{3}(\phi^\star \nabla \mu_\phi + \frac{2\delta t}{3}\rho^\star \nabla \mu_\rho, \boldsymbol{v}) = (\boldsymbol{f}_1, \boldsymbol{v}), \tag{3.38}$$

$$(\phi, w) - \frac{2\delta t}{3}(\boldsymbol{u}\phi^\star, \nabla w) + \frac{2\delta t}{3}M_1(\nabla \mu, \nabla w) = (f_2, w), \tag{3.39}$$

$$(-\mu_\phi, \psi) + \epsilon(\nabla \phi, \nabla \psi) + \frac{1}{\epsilon}2(\phi^\star \phi, \phi^\star \psi) - \alpha(\rho \boldsymbol{Z}^\star, \nabla \psi) + \alpha(\boldsymbol{Z}^\star \cdot \nabla \phi, \boldsymbol{Z}^\star \cdot \nabla \psi) = (f_3, \psi), \tag{3.40}$$

$$(\rho, \varpi) - \frac{2\delta t}{3}(\boldsymbol{u}\rho^\star, \nabla \varpi) + \frac{2\delta t}{3}M_2(\nabla \mu_\rho, \nabla \varpi) = (f_4, \varpi), \tag{3.41}$$

$$(-\mu_\rho, \varrho) + \eta(\nabla \rho, \nabla \varrho) + \alpha(\rho, \varrho) - \alpha(\boldsymbol{Z}^\star \cdot \nabla \phi, \varrho) + \frac{1}{2}\beta(H^\star \rho, H^\star \varrho) = (f_5, \varrho), \tag{3.42}$$

for any $\boldsymbol{v} \in H^1_{\boldsymbol{u}}(\Omega)$ and $w, \psi, \varpi, \varrho \in \bar{H}^1(\Omega)$, where $\boldsymbol{f}_1$ and $f_{i,i=2,\cdots,5}$ are the corresponding term on the right hand side of (3.31)-(3.35).

We denote the above linear system (3.38)-(3.42) as

$$(\mathbb{L}(\boldsymbol{X}), \boldsymbol{Y}) = (\mathbb{B}, \boldsymbol{Y}) \tag{3.43}$$

where $\mathbb{L}$ is the linear operator, $\boldsymbol{X} = (\boldsymbol{u}, \mu_\phi, \phi, \mu_\rho, \rho)^T$, and $\boldsymbol{Y} = (\boldsymbol{v}, w, \psi, \varpi, \varrho)^T$, and $\boldsymbol{X}, \boldsymbol{Y} \in (H^1_{\boldsymbol{u}}, \bar{H}^1, \bar{H}^1, \bar{H}^1, \bar{H}^1)(\Omega)$.

**Theorem 3.1.** *The linear system* (3.43) *admits a unique solution* $\boldsymbol{X} = (\boldsymbol{u}, \mu_\phi, \phi, \mu_\rho, \rho)^T \in (H^1_{\boldsymbol{u}}, \bar{H}^1, \bar{H}^1, \bar{H}^1, \bar{H}^1)(\Omega)$.

*Proof.* (i). For any $\boldsymbol{X} = (\boldsymbol{u}, \mu_\phi, \phi, \mu_\rho, \rho)^T$ and $\boldsymbol{Y} = (\boldsymbol{v}, w, \psi, \varpi, \varrho)^T$ with $\boldsymbol{X}, \boldsymbol{Y} \in (H^1_{\boldsymbol{u}}, \bar{H}^1, \bar{H}^1, \bar{H}^1, \bar{H}^1)(\Omega)$, we have

$$(\mathbb{L}(\boldsymbol{X}), \boldsymbol{Y}) \le C_1(\|\boldsymbol{u}\|_{H^1} + \|\mu_\phi\|_{H^1} + \|\phi\|_{H^1} + \|\mu_\rho\|_{H^1} + \|\rho\|_{H^1})$$
$$(\|\boldsymbol{v}\|_{H^1} + \|w\|_{H^1} + \|\psi\|_{H^1} + \|\varpi\|_{H^1} + \|\varrho\|_{H^1}), \tag{3.44}$$

where $C_1$ is a constant that depends on $\delta t, \nu, M_1, M_2, \epsilon, \alpha, \beta, \eta$, $\|\phi^\star\|_{L^\infty}$, $\|\boldsymbol{Z}^\star\|_{L^\infty}$, $\|\rho^\star\|_{L^\infty}$, and $\|H^\star\|_{L^\infty}$.

(ii). It is easy to derive that

$$
\begin{aligned}
(\mathbb{L}(\boldsymbol{X}),\boldsymbol{X}) &= \|\boldsymbol{u}\|^2 + \frac{2\delta t}{3}\nu\|\nabla\boldsymbol{u}\|^2 + \frac{2\delta t}{3}M_1\|\nabla\mu_\phi\|^2 + \frac{2\delta t}{3}M_2\|\nabla\mu_\rho\|^2\\
&\quad + \epsilon\|\nabla\phi\|^2 + \eta\|\nabla\rho\|^2 + \frac{2}{\epsilon}\|\phi^\star\phi\|^2 + \frac{\beta}{2}\|H^\star\rho\|^2 + \alpha\|\rho - \boldsymbol{Z}^\star\cdot\nabla\phi\|^2\\
&\geq C_2(\|\boldsymbol{u}\|^2_{H^1} + \|\phi\|^2_{H^1} + \|\rho\|^2_{H^1} + \|\mu_\phi\|^2_{H^1} + \|\mu_\rho\|^2_{H^1}),
\end{aligned} \tag{3.45}
$$

where $C_2$ is a constant that depends on $\delta t, \nu, M_1, M_2, \epsilon, \eta$.

Then from the Lax-Milgram theorem, we conclude that the linear system (3.43) admits a unique solution $\boldsymbol{X} = (\boldsymbol{u}, \mu_\phi, \phi, \mu_\rho, \rho)^T \in (H^1_{\boldsymbol{u}}, \bar{H}^1, \bar{H}^1, \bar{H}^1, \bar{H}^1)(\Omega)$. □

The stability result of the BDF2 scheme follows the same lines as in the derivation of the new PDE energy dissipation law (3.17), as follows.

**Theorem 3.2.** *The second order linear scheme* (3.18)–(3.28) *is unconditionally energy stable, i.e., it satisfies the following discrete energy dissipation law:*

$$
\frac{1}{\delta t}(E^{n+1,n}_{2nd} - E^{n,n-1}_{2nd}) \leq -\nu\|\nabla\tilde{\boldsymbol{u}}^{n+1}\|^2 - M_1\|\nabla\mu^{n+1}_\phi\|^2 - M_2\|\nabla\mu^{n+1}_\rho\|^2, \tag{3.46}
$$

*where*

$$
\begin{aligned}
E^{n+1,n}_{2nd} &= \frac{1}{2}\Big(\frac{1}{2}\|\boldsymbol{u}^{n+1}\|^2 + \frac{1}{2}\|2\boldsymbol{u}^{n+1} - \boldsymbol{u}^n\|^2\Big) + \frac{\delta t^2}{3}\|\nabla p^{n+1}\|^2\\
&\quad + \frac{\epsilon}{2}\Big(\frac{1}{2}\|\nabla\phi^{n+1}\|^2 + \frac{1}{2}\|2\nabla\phi^{n+1} - \nabla\phi^n\|^2\Big) + \frac{\eta}{2}\Big(\frac{1}{2}\|\nabla\rho^{n+1}\|^2 + \frac{1}{2}\|2\nabla\rho^{n+1} - \nabla\rho^n\|^2\Big)\\
&\quad + \frac{1}{4\epsilon}\Big(\frac{1}{2}\|U^{n+1}\|^2 + \frac{1}{2}\|2U^{n+1} - U^n\|^2\Big) + \frac{\alpha}{2}\Big(\frac{1}{2}\|V^{n+1}\|^2 + \frac{1}{2}\|2V^{n+1} - V^n\|^2\Big)\\
&\quad + \beta\Big(\frac{1}{2}\|W^{n+1}\|^2 + \frac{1}{2}\|2W^{n+1} - W^n\|^2\Big).
\end{aligned} \tag{3.47}
$$

*Proof.* By taking the $L^2$ inner product of (3.18) with $2\delta t\tilde{\boldsymbol{u}}^{n+1}$, we obtain

$$
\begin{aligned}
&(3\tilde{\boldsymbol{u}}^{n+1} - 4\boldsymbol{u}^n + \boldsymbol{u}^{n-1}, \tilde{\boldsymbol{u}}^{n+1}) + 2\nu\delta t\|\nabla\tilde{\boldsymbol{u}}^{n+1}\|^2 + 2\delta t(\nabla p^n, \tilde{\boldsymbol{u}}^{n+1})\\
&\qquad + 2\delta t(\phi^\star\nabla\mu^{n+1}_\phi, \tilde{\boldsymbol{u}}^{n+1}) + 2\delta t(\rho^\star\nabla\mu^{n+1}_\rho, \tilde{\boldsymbol{u}}^{n+1}) = 0.
\end{aligned} \tag{3.48}
$$

From (3.27), for any function $\boldsymbol{v}$ with $\nabla\cdot\boldsymbol{v} = 0$, we can derive

$$
(\boldsymbol{u}^{n+1}, \boldsymbol{v}) = (\tilde{\boldsymbol{u}}^{n+1}, \boldsymbol{v}). \tag{3.49}
$$

Then for the first term in (3.48), we have

$$
\begin{aligned}
&(3\tilde{\boldsymbol{u}}^{n+1} - 4\boldsymbol{u}^n + \boldsymbol{u}^{n-1}, \tilde{\boldsymbol{u}}^{n+1})\\
&= (3\tilde{\boldsymbol{u}}^{n+1} - 4\boldsymbol{u}^n + \boldsymbol{u}^{n-1}, \boldsymbol{u}^{n+1}) + (3\tilde{\boldsymbol{u}}^{n+1} - 4\boldsymbol{u}^n + \boldsymbol{u}^{n-1}, \tilde{\boldsymbol{u}}^{n+1} - \boldsymbol{u}^{n+1})\\
&= (3\boldsymbol{u}^{n+1} - 4\boldsymbol{u}^n + \boldsymbol{u}^{n-1}, \boldsymbol{u}^{n+1}) + (3\tilde{\boldsymbol{u}}^{n+1}, \tilde{\boldsymbol{u}}^{n+1} - \boldsymbol{u}^{n+1})\\
&= (3\boldsymbol{u}^{n+1} - 4\boldsymbol{u}^n + \boldsymbol{u}^{n-1}, \boldsymbol{u}^{n+1}) + 3(\tilde{\boldsymbol{u}}^{n+1} - \boldsymbol{u}^{n+1}, \tilde{\boldsymbol{u}}^{n+1} + \boldsymbol{u}^{n+1})\\
&= \frac{1}{2}\Big(\|\boldsymbol{u}^{n+1}\|^2 - \|\boldsymbol{u}^n\|^2 + \|2\boldsymbol{u}^{n+1} - \boldsymbol{u}^n\|^2 - \|2\boldsymbol{u}^n - \boldsymbol{u}^{n-1}\|^2 + \|\boldsymbol{u}^{n+1} - 2\boldsymbol{u}^n + \boldsymbol{u}^{n-1}\|^2\Big)\\
&\quad + 3(\|\tilde{\boldsymbol{u}}^{n+1}\|^2 - \|\boldsymbol{u}^{n+1}\|^2),
\end{aligned} \tag{3.50}
$$

where we have use the following identity

$$2(3a-4b+c,a)=|a|^2-|b|^2+|2a-b|^2-|2b-c|^2+|a-2b+c|^2. \tag{3.51}$$

For the projection step, we rewrite (3.27) as

$$\frac{3}{2\delta t}\boldsymbol{u}^{n+1}+\nabla p^{n+1}=\frac{3}{2\delta t}\tilde{\boldsymbol{u}}^{n+1}+\nabla p^{n}. \tag{3.52}$$

By squaring both sides of the above equality, we obtain

$$\frac{9}{4\delta t^2}\|\boldsymbol{u}^{n+1}\|^2+\|\nabla p^{n+1}\|^2=\frac{9}{4\delta t^2}\|\tilde{\boldsymbol{u}}^{n+1}\|^2+\|\nabla p^{n}\|^2+\frac{3}{\delta t}(\tilde{\boldsymbol{u}}^{n+1},\nabla p^n), \tag{3.53}$$

namely, by multiplying $2\delta t^2/3$, we have

$$\frac{3}{2}(\|\boldsymbol{u}^{n+1}\|^2-\|\tilde{\boldsymbol{u}}^{n+1}\|^2)+\frac{2\delta t^2}{3}(\|\nabla p^{n+1}\|^2-\|\nabla p^{n}\|^2)=2\delta t(\tilde{\boldsymbol{u}}^{n+1},\nabla p^n). \tag{3.54}$$

By taking the $L^2$ inner product of (3.27) with $2\delta t\boldsymbol{u}^{n+1}$, we have

$$\frac{3}{2}(\|\boldsymbol{u}^{n+1}\|^2-\|\tilde{\boldsymbol{u}}^{n+1}\|^2+\|\boldsymbol{u}^{n+1}-\tilde{\boldsymbol{u}}^{n+1}\|^2)=0. \tag{3.55}$$

By combining (3.48), (3.50), (3.54), and (3.55), we obtain

$$\begin{aligned}&\frac{1}{2}(\|\boldsymbol{u}^{n+1}\|^2-\|\boldsymbol{u}^{n}\|^2+\|2\boldsymbol{u}^{n+1}-\boldsymbol{u}^n\|^2-\|2\boldsymbol{u}^{n}-\boldsymbol{u}^{n-1}\|^2+\|\boldsymbol{u}^{n+1}-2\boldsymbol{u}^n+\boldsymbol{u}^{n-1}\|^2)\\&\quad+\frac{3}{2}\|\boldsymbol{u}^{n+1}-\tilde{\boldsymbol{u}}^{n+1}\|^2+\frac{2\delta t^2}{3}(\|\nabla p^{n+1}\|^2-\|\nabla p^n\|^2)+2\nu\delta t\|\nabla\tilde{\boldsymbol{u}}^{n+1}\|^2\\&\quad+2\delta t(\phi^\star\nabla\mu_\phi^{n+1},\tilde{\boldsymbol{u}}^{n+1})+2\delta t(\rho^\star\nabla\mu_\rho^{n+1},\tilde{\boldsymbol{u}}^{n+1})=0.\end{aligned} \tag{3.56}$$

By taking the $L^2$ inner product of (3.19) with $2\delta t\mu_\phi^{n+1}$, we obtain

$$(3\phi^{n+1}-4\phi^n+\phi^{n-1},\mu_\phi^{n+1})-2\delta t(\phi^\star\tilde{\boldsymbol{u}}^{n+1},\nabla\mu_\phi^{n+1})=-2M_1\delta t\|\nabla\mu_\phi^{n+1}\|^2. \tag{3.57}$$

By taking the $L^2$ inner product of (3.20) with $-(3\phi^{n+1}-4\phi^n+\phi^{n-1})$ and applying (3.51), we obtain

$$\begin{aligned}&-(\mu_\phi^{n+1},3\phi^{n+1}-4\phi^n+\phi^{n-1})\\&\quad=-\frac{\epsilon}{2}\Big(\|\nabla\phi^{n+1}\|^2-\|\nabla\phi^n\|^2+\|2\nabla\phi^{n+1}-\nabla\phi^n\|^2-\|2\nabla\phi^n-\nabla\phi^{n-1}\|^2\\&\quad+\|\nabla\phi^{n+1}-2\nabla\phi^n+\nabla\phi^{n-1}\|^2\Big)-\frac{1}{\epsilon}(\phi^\star U^{n+1},3\phi^{n+1}-4\phi^n+\phi^{n-1})\\&\quad+\alpha(V^{n+1}\boldsymbol{Z}^\star,\nabla(3\phi^{n+1}-4\phi^n+\phi^{n-1})).\end{aligned} \tag{3.58}$$

By taking the $L^2$ inner product of (3.21) with $2\delta t\mu_\rho^{n+1}$, we obtain

$$(3\rho^{n+1}-4\rho^n+\rho^{n-1},\mu_\rho^{n+1})-2\delta t(\rho^\star\tilde{\boldsymbol{u}}^{n+1},\nabla\mu_\rho^{n+1})=-2M_2\delta t\|\nabla\mu_\rho^{n+1}\|^2. \tag{3.59}$$

By taking the $L^2$ inner product of (3.22) with $-(3\rho^{n+1}-4\rho^n+\rho^{n-1})$, we obtain

$$
\begin{aligned}
&-(\mu_\rho^{n+1},3\rho^{n+1}-4\rho^n+\rho^{n-1})\\
&\quad=-\frac{\eta}{2}\Big(\|\nabla\rho^{n+1}\|^2-\|\nabla\rho^n\|^2+\|2\nabla\rho^{n+1}-\nabla\rho^n\|^2-\|2\nabla\rho^n-\nabla\rho^{n-1}\|^2\\
&\quad\quad+\|\nabla\rho^{n+1}-2\nabla\rho^n+\nabla\rho^{n-1}\|^2\Big)-\alpha(V^{n+1},3\rho^{n+1}-4\rho^n+\rho^{n-1})\\
&\quad\quad-\beta(H^\star W^{n+1},3\rho^{n+1}-4\rho^n+\rho^{n-1}).
\end{aligned}
\tag{3.60}
$$

By taking the $L^2$ inner product of (3.23) with $\frac{1}{2\epsilon}U^{n+1}$, we obtain

$$
\begin{aligned}
&\frac{1}{4\epsilon}\Big(\|U^{n+1}\|^2-\|U^n\|^2+\|2U^{n+1}-U^n\|^2-\|2U^n-U^{n-1}\|^2+\|U^{n+1}-2U^n+U^{n-1}\|^2\Big)\\
&\quad=\frac{1}{\epsilon}(\phi^\star(3\phi^{n+1}-4\phi^n+\phi^{n-1}),U^{n+1}).
\end{aligned}
\tag{3.61}
$$

By taking the $L^2$ inner product of (3.24) with $\alpha V^{n+1}$, we obtain

$$
\begin{aligned}
&\frac{\alpha}{2}(\|V^{n+1}\|^2-\|V^n\|^2+\|2V^{n+1}-V^n\|^2-\|2V^n-V^{n-1}\|^2+\|V^{n+1}-2V^n+V^{n-1}\|^2)\\
&\quad=\alpha(3\rho^{n+1}-4\rho^n+\rho^{n-1},V^{n+1})-\alpha(\boldsymbol{Z}^\star\nabla(3\phi^{n+1}-4\phi^n+\phi^{n-1}),V^{n+1}).
\end{aligned}
\tag{3.62}
$$

By taking the $L^2$ inner product of (3.25) with $2\beta W^{n+1}$, we obtain

$$
\begin{aligned}
&\beta(\|W^{n+1}\|^2-\|W^n\|^2+\|2W^{n+1}-W^n\|^2-\|2W^n-W^{n-1}\|^2\\
&\quad+\|W^{n+1}-2W^n+W^{n-1}\|^2)=\beta(H^\star(3\rho^{n+1}-4\rho^n+\rho^{n-1}),W^{n+1}).
\end{aligned}
\tag{3.63}
$$

Combination of (3.56), (3.57)–(3.63) gives us

$$
\begin{aligned}
&\frac{1}{2}(\|\boldsymbol{u}^{n+1}\|^2-\|\boldsymbol{u}^n\|^2+\|2\boldsymbol{u}^{n+1}-\boldsymbol{u}^n\|^2-\|2\boldsymbol{u}^n-\boldsymbol{u}^{n-1}\|^2+\|\boldsymbol{u}^{n+1}-2\boldsymbol{u}^n+\boldsymbol{u}^{n-1}\|^2)\\
&\quad+\frac{3}{2}\|\boldsymbol{u}^{n+1}-\tilde{\boldsymbol{u}}^{n+1}\|^2+\frac{2\delta t^2}{3}(\|\nabla p^{n+1}\|^2-\|\nabla p^n\|^2)\\
&\quad+\frac{\epsilon}{2}\Big(\|\nabla\phi^{n+1}\|^2-\|\nabla\phi^n\|^2+\|2\nabla\phi^{n+1}-\nabla\phi^n\|^2-\|2\nabla\phi^n-\nabla\phi^{n-1}\|^2\\
&\quad+\|\nabla\phi^{n+1}-2\nabla\phi^n+\nabla\phi^{n-1}\|^2\Big)\\
&\quad+\frac{\eta}{2}\Big(\|\nabla\rho^{n+1}\|^2-\|\nabla\rho^n\|^2+\|2\nabla\rho^{n+1}-\nabla\rho^n\|^2-\|2\nabla\rho^n-\nabla\rho^{n-1}\|^2\\
&\quad+\|\nabla\rho^{n+1}-2\nabla\rho^n+\nabla\rho^{n-1}\|^2\Big)\\
&\quad+\frac{1}{4\epsilon}\Big(\|U^{n+1}\|^2-\|U^n\|^2+\|2U^{n+1}-U^n\|^2-\|2U^n-U^{n-1}\|^2+\|U^{n+1}-2U^n+U^{n-1}\|^2\Big)\\
&\quad+\frac{\alpha}{2}\Big(\|V^{n+1}\|^2-\|V^n\|^2+\|2V^{n+1}-V^n\|^2-\|2V^n-V^{n-1}\|^2+\|V^{n+1}-2V^n+V^{n-1}\|^2\Big)\\
&\quad+\beta\Big(\|W^{n+1}\|^2-\|W^n\|^2+\|2W^{n+1}-W^n\|^2-\|2W^n-W^{n-1}\|^2+\|W^{n+1}-2W^n+W^{n-1}\|^2\Big)\\
&=-2\nu\delta t\|\nabla\tilde{\boldsymbol{u}}^{n+1}\|^2-2M_1\delta t\|\nabla\mu_\phi^{n+1}\|^2-2M_2\delta t\|\nabla\mu_\rho^{n+1}\|^2.
\end{aligned}
$$

Finally, we obtain the result (3.46) after dropping some positive terms from the above equation. □

3.2. **Crank-Nicolson Scheme.** One also can easily develop an alternative version of second order scheme based on the Crank-Nicolson type scheme. It reads as follows.

Assuming that $(\boldsymbol{u},\phi,\rho,U,V,W)^{n-1}$ and $(\boldsymbol{u},\phi,\rho,U,V,W)^{n}$ are known,

**step 1:** we update $\tilde{\boldsymbol{u}}^{n+1}$, $\phi^{n+1}$, $\rho^{n+1}$, $U^{n+1}$, $V^{n+1}$, $W^{n+1}$ from

$$\frac{\tilde{\boldsymbol{u}}^{n+1}-\boldsymbol{u}^n}{\delta t}+B(\boldsymbol{u}^\flat,\tilde{\boldsymbol{u}}^{n+\frac12})-\nu\Delta\tilde{\boldsymbol{u}}^{n+\frac12}+\nabla p^n+\phi^\flat\nabla\mu_\phi^{n+\frac12}+\rho^\flat\nabla\mu_\rho^{n+\frac12}=0 \tag{3.64}$$

$$\frac{\phi^{n+1}-\phi^n}{\delta t}+\nabla\cdot(\tilde{\boldsymbol{u}}^{n+\frac12}\phi^\flat)=M_1\Delta\mu_\phi^{n+\frac12}, \tag{3.65}$$

$$\mu_\phi^{n+\frac12}=-\epsilon\Delta\phi^{n+\frac12}+\frac{1}{\epsilon}\phi^\flat U^{n+\frac12}+\alpha\nabla\cdot(V^{n+\frac12}\boldsymbol{Z}^\flat), \tag{3.66}$$

$$\frac{\rho^{n+1}-\rho^n}{\delta t}+\nabla\cdot(\rho^\flat\tilde{\boldsymbol{u}}^{n+\frac12})=M_2\Delta\mu_\rho^{n+\frac12}, \tag{3.67}$$

$$\mu_\rho^{n+\frac12}=-\eta\Delta\rho^{n+\frac12}+\alpha V^{n+\frac12}+\beta H^\flat W^{n+\frac12}, \tag{3.68}$$

$$U^{n+1}-U^n=2\phi^\flat(\phi^{n+1}-\phi^n), \tag{3.69}$$

$$V^{n+1}-V^n=(\rho^{n+1}-\rho^n)-\boldsymbol{Z}^\flat\cdot\nabla(\phi^{n+1}-\phi^n), \tag{3.70}$$

$$W^{n+1}-W^n=\frac12 H^\flat(\rho^{n+1}-\rho^n), \tag{3.71}$$

with periodic boundary condition being imposed, where

$$\begin{aligned}
&\tilde{\boldsymbol{u}}^{n+\frac12}=\frac{\tilde{\boldsymbol{u}}^{n+1}+\boldsymbol{u}^n}{2},\\
&\phi^{n+\frac12}=\frac{\phi^{n+1}+\phi^n}{2},\rho^{n+\frac12}=\frac{\rho^{n+1}+\rho^n}{2},\\
&U^{n+\frac12}=\frac{U^{n+1}+U^n}{2},V^{n+\frac12}=\frac{V^{n+1}+V^n}{2},W^{n+\frac12}=\frac{W^{n+1}+W^n}{2},\\
&\phi^\flat=\frac32\phi^n-\frac12\phi^{n-1},\rho^\flat=\frac32\rho^n-\frac12\rho^{n-1},\boldsymbol{u}^\flat=\frac32\boldsymbol{u}^n-\frac12\boldsymbol{u}^{n-1},\\
&\boldsymbol{Z}^\flat=\boldsymbol{Z}(\frac32\phi^n-\frac12\phi^{n-1}),H^\flat=H(\frac32\rho^n-\frac12\rho^{n-1}).
\end{aligned} \tag{3.72}$$

.

**step 2:** we update $\boldsymbol{u}^{n+1}$ and $p^{n+1}$ from

$$\frac{1}{\delta t}(\boldsymbol{u}^{n+1}-\tilde{\boldsymbol{u}}^{n+1})+\frac12(\nabla p^{n+1}-\nabla p^n)=0, \tag{3.73}$$

$$\nabla\cdot\boldsymbol{u}^{n+1}=0, \tag{3.74}$$

with the periodic boundary conditions.

The well-posedness and unconditionally energy stability of the Crank-Nicolson scheme can be easily proved in the similar way of handling the BDF2 scheme. Therefore we omit the details and leave them to the interested readers.

## 4. Numerical examples

We now present numerical experiments in two and three dimensions to validate the theoretical results derived in the previous section and demonstrate the efficiency, energy stability and accuracy of the proposed numerical schemes. In all examples, we set the domain $\Omega=[0,2\pi]^d, d=2,3$. If not explicitly

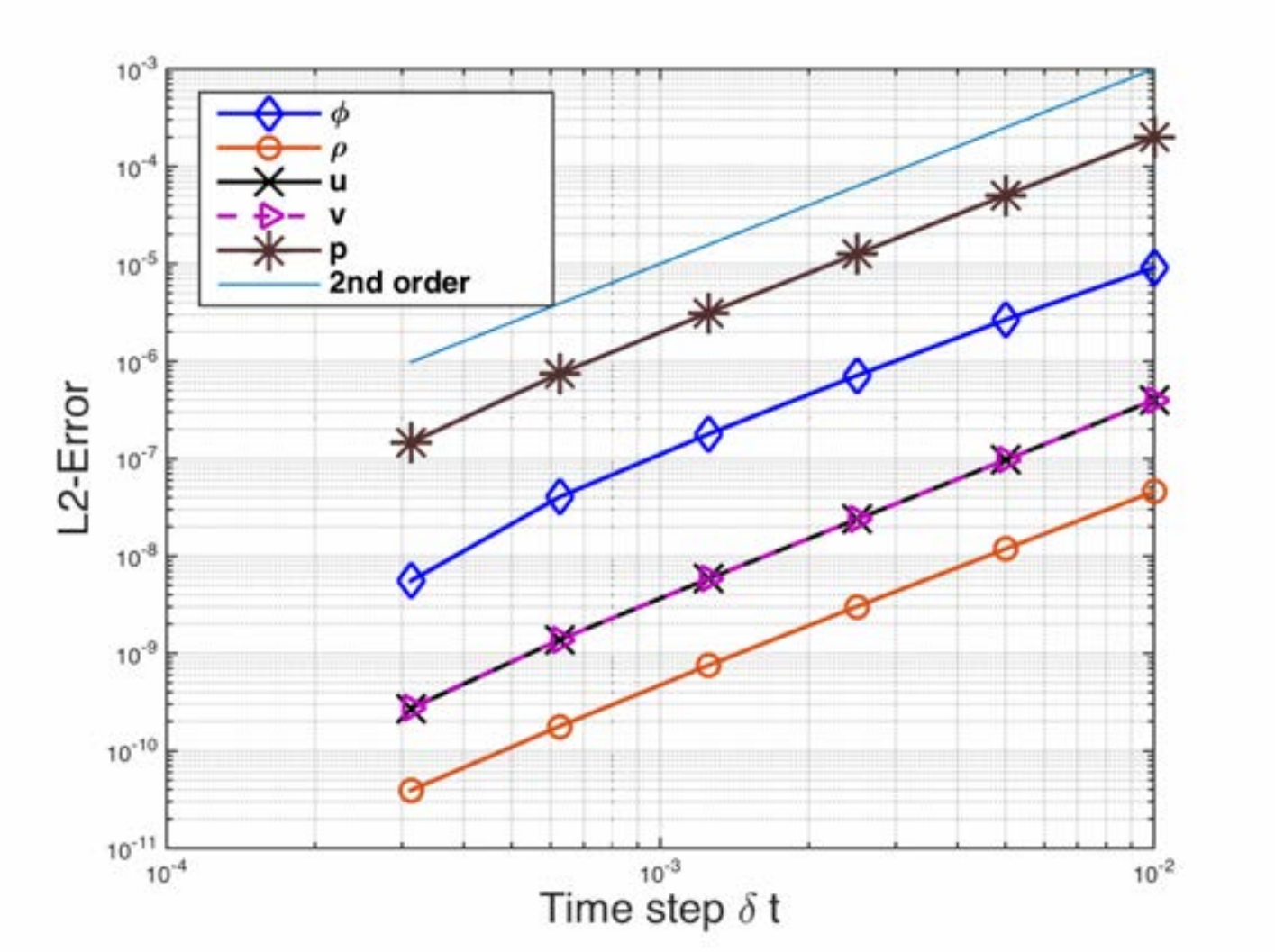

FIGURE. 4.1. The $L^2$ numerical errors of all variables $\phi, \rho, \boldsymbol{u} = (u, v), p$ at $t = 0.2$ that are computed using the scheme BDF2 and various temporal resolutions with the initial conditions of (4.2) for mesh refinement test in time.

specified, the default values of order parameters are given as follows,

$$\widehat{\epsilon} = 1\text{e}{-4},\ M_\mu = M_\rho = 1\text{e}{-2},\ \epsilon = 5\text{e}{-2},\ \alpha = 1\text{e}{-2},\ \beta = 5\text{e}{-2},\ \nu = 1,\ A = 1. \tag{4.1}$$

We use the Fourier-spectral method to discretize the space, and $129^d$ Fourier modes are used so that the errors from the spatial discretization is negligible compared with the time discretization errors.

4.1. **Accuracy test.** We first test convergence rates of proposed BDF2 scheme (3.18)-(3.28). The following initial conditions

$$\begin{cases} \phi_0(x, y) = 0.1\cos(3x) + 0.4\cos(y), \\ \rho_0(x, y) = 0.2\sin(2x) + 0.5\sin(y) \end{cases} \tag{4.2}$$

are used. We perform the refinement test of the time step size, and choose the approximate solution obtained by using the scheme with the time step size $\delta t = 1.5625\text{e}{-4}$ as the benchmark solution (approximately the exact solution) for computing errors. We present the $L^2$ errors of all variables between the numerical solution and the exact solution at $t = 0.2$ with different time step sizes in Fig. 4.1. We observe that the scheme are almost perfect second order accurate in time.

4.2. **Coarsening effect of two circles.** In this example, we test the scheme BDF2 by using the initial conditions of two circles with different radii to see how the coarsening effects execute. We set the initial condition as

$$\begin{cases} \phi_0(x, y) = \sum_{i=1}^{2} -\tanh\Big(\dfrac{\sqrt{(x - x_i)^2 + (y - y_i)^2 - r_i}}{1.2\epsilon}\Big) + 1, \\ \rho_0(x, y) = 0.3, \end{cases} \tag{4.3}$$

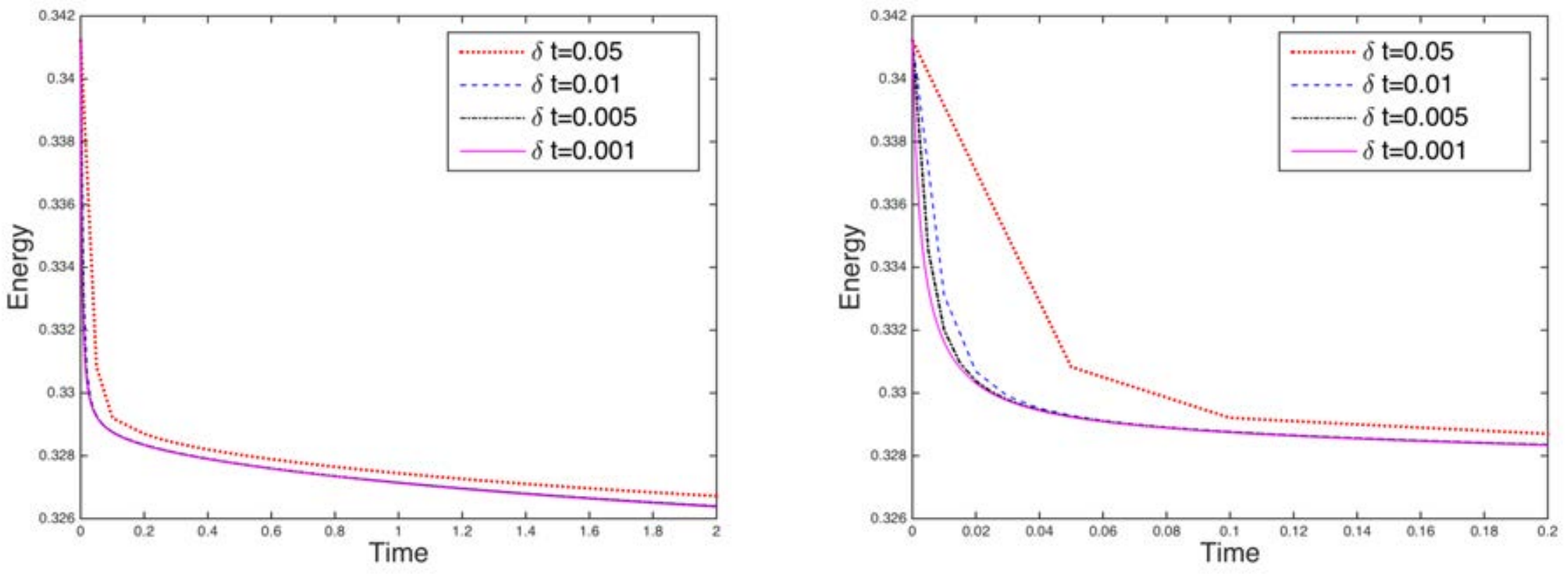

Figure. 4.2. Time evolution of the free energy functional (3.47) for four different time step sizes, $\delta t$ = 5e−2, 1e−2, 5e−3, and 1e−3 for the Example 4.2 with default parameters (4.1). The energy curves show the decays for all time step sizes, that confirms that our algorithm is unconditionally stable. The left subfigure is the energy profile for $t \in [0, 2]$, and the right subfigure is a close-up view for $t \in [0, 0.2]$.

where $(x_1, y_1, r_1) = (\pi - 0.7, \pi - 0.6, 1.5)$ and $(x_2, y_2, r_2) = (\pi + 1.65, \pi + 1.6, 0.7)$.

We emphasize that any time step size $\delta t$ is allowable for the computations from the stability concern since all developed schemes are unconditionally energy stable. But larger time step will definitely induce large numerical errors. Therefore, we need to discover the rough range of the allowable maximum time step size in order to obtain good accuracy and to consume as low computational cost as possible. This time step range could be estimated through the energy evolution curve plots, shown in Fig. 4.2, where we compare the time evolution of the discrete free energy (3.47) for four different time step sizes until $t = 2$ using the scheme BDF2. We observe that all four energy curves decay monotonically, which numerically confirms that our algorithms are unconditionally energy stable. For smaller time steps of $\delta t$ = 1e−3, 5e−2 and 1e−2, the three energy curves coincide very well. But for the largest time step of $\delta t$ = 5e−2, the energy curve deviates viewable away from others. This means the adopted time step size should not be larger than 5e−2, in order to get reasonably good accuracy.

In Fig. 4.3, we show the evolutions of the phase field variable $\phi$ and concentration variable $\rho$ at various time by using the time step $\delta t$ = 1e−3. we observe the coarsening effect that the small circle is absorbed into the big circle, and the total absorption happens at around $t = 98$. The velocity field $\boldsymbol{u}$, together with the interface contours of $\phi$, are plotted in Fig. 4.4. The time evolutions of the two free energy functionals, the modified free energy (3.47) (in terms of the three new variables) and the original energy (3.2), are plotted in Fig. 4.5. It is remarkable that these two energies actually coincide very well, and both of them decay to the equilibrium, monotonically. At around $t = 100$, the energies undergo a rapid decrease when the small circle is compeletly absorbed. Soon after that, the system achieves the equilibrium of circular shape.

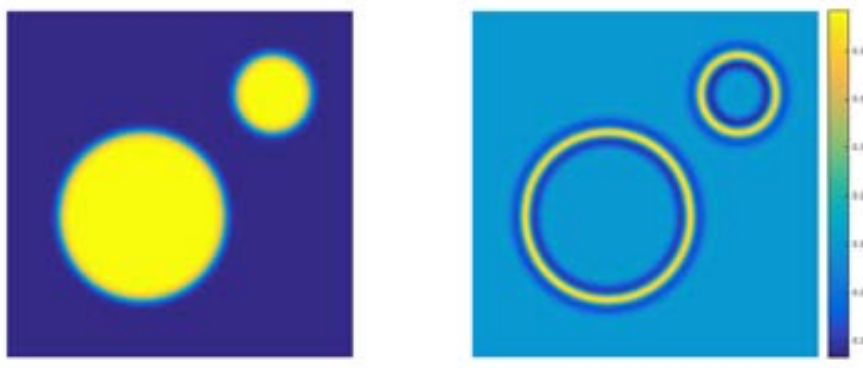

(a) $t = 1$.

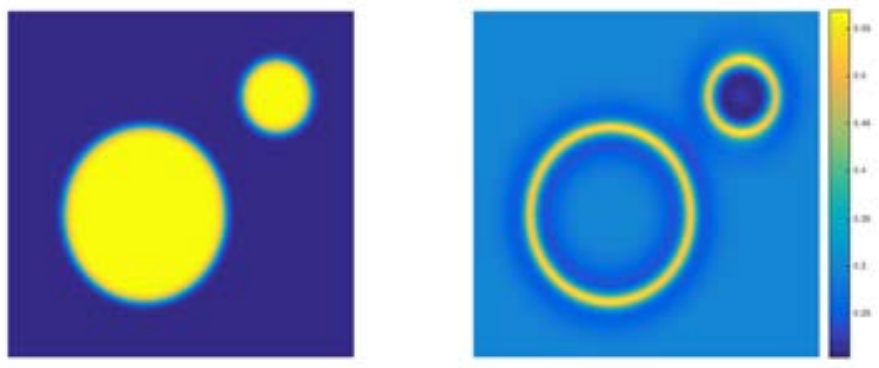

(b) $t = 20$.

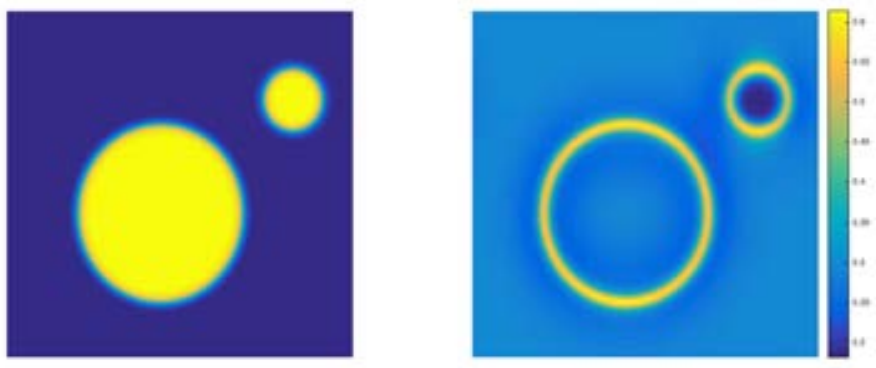

(c) $t = 50$.

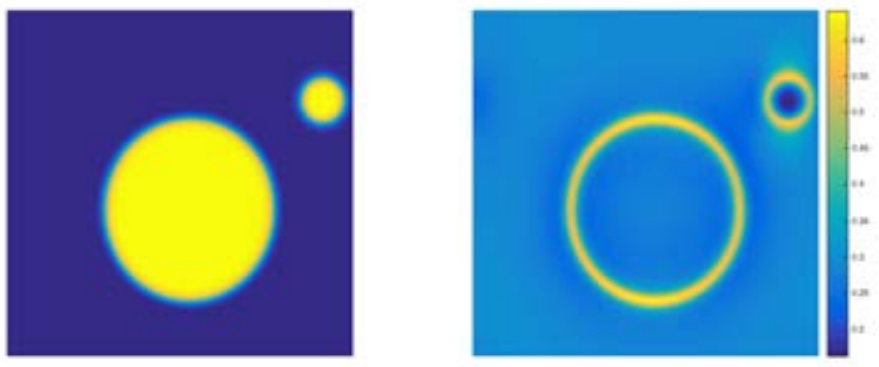

(d) $t = 80$.

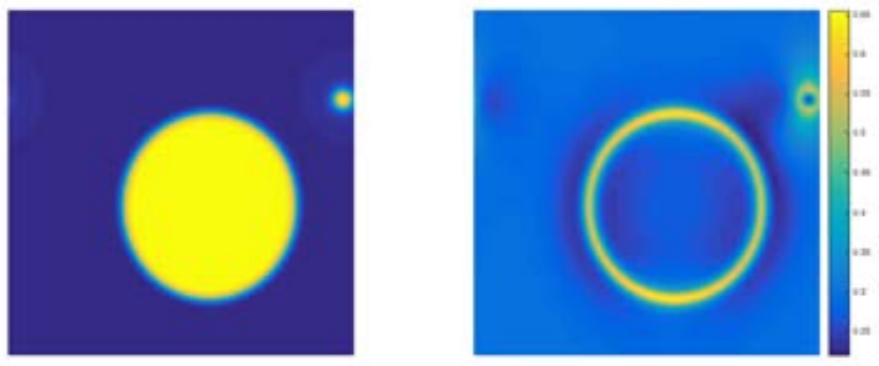

(e) $t = 98$.

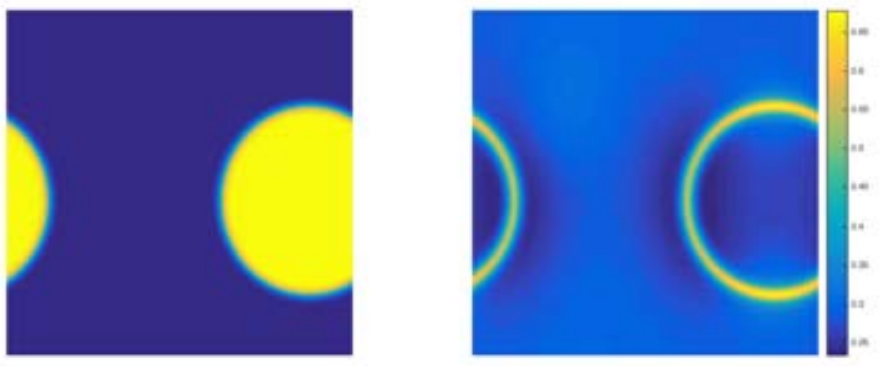

(f) $t = 200$.

FIGURE. 4.3. Snapshots of the phase variables $\phi$ and $\rho$ are taken at $t = 1, 20, 50, 80, 98, 200$ for Example 4.2 by using the initial condition (4.3). For each panel, the left subfigure is the profile of $\phi$ and the right subfigure is the profile of $\rho$.

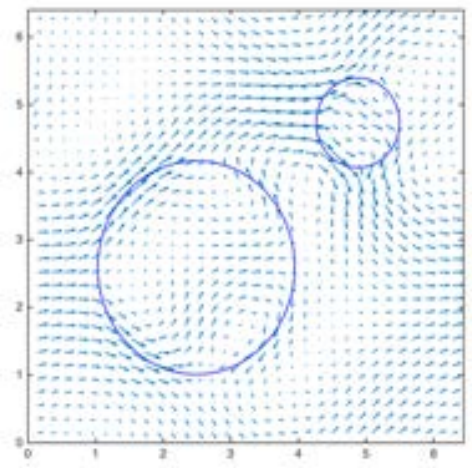

(a) $t = 20$.

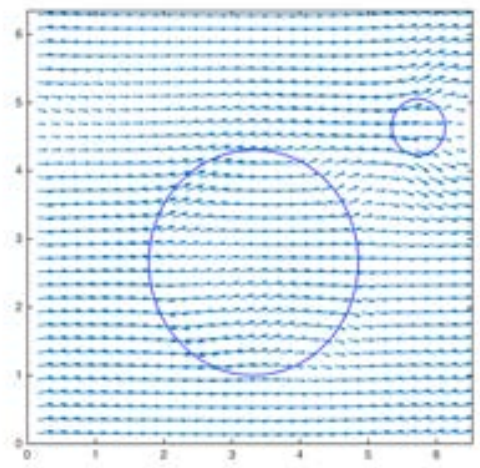

(b) $t = 80$.

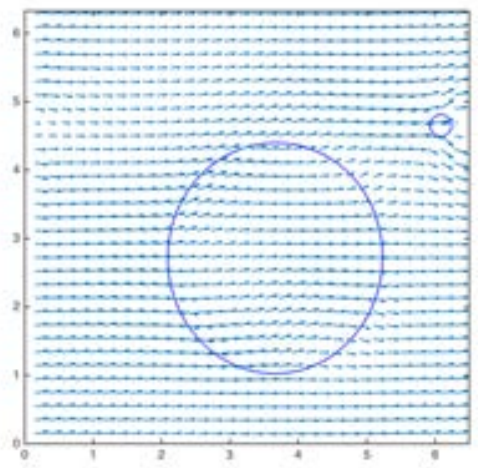

(c) $t = 98$.

FIGURE. 4.4. Profiles of the velocity field $\boldsymbol{u}$ together with the interface contour $\{\phi = 0\}$ for Example 4.2 by using the initial condition (4.3). Snapshots are taken at $t = 20$, 80, 98.

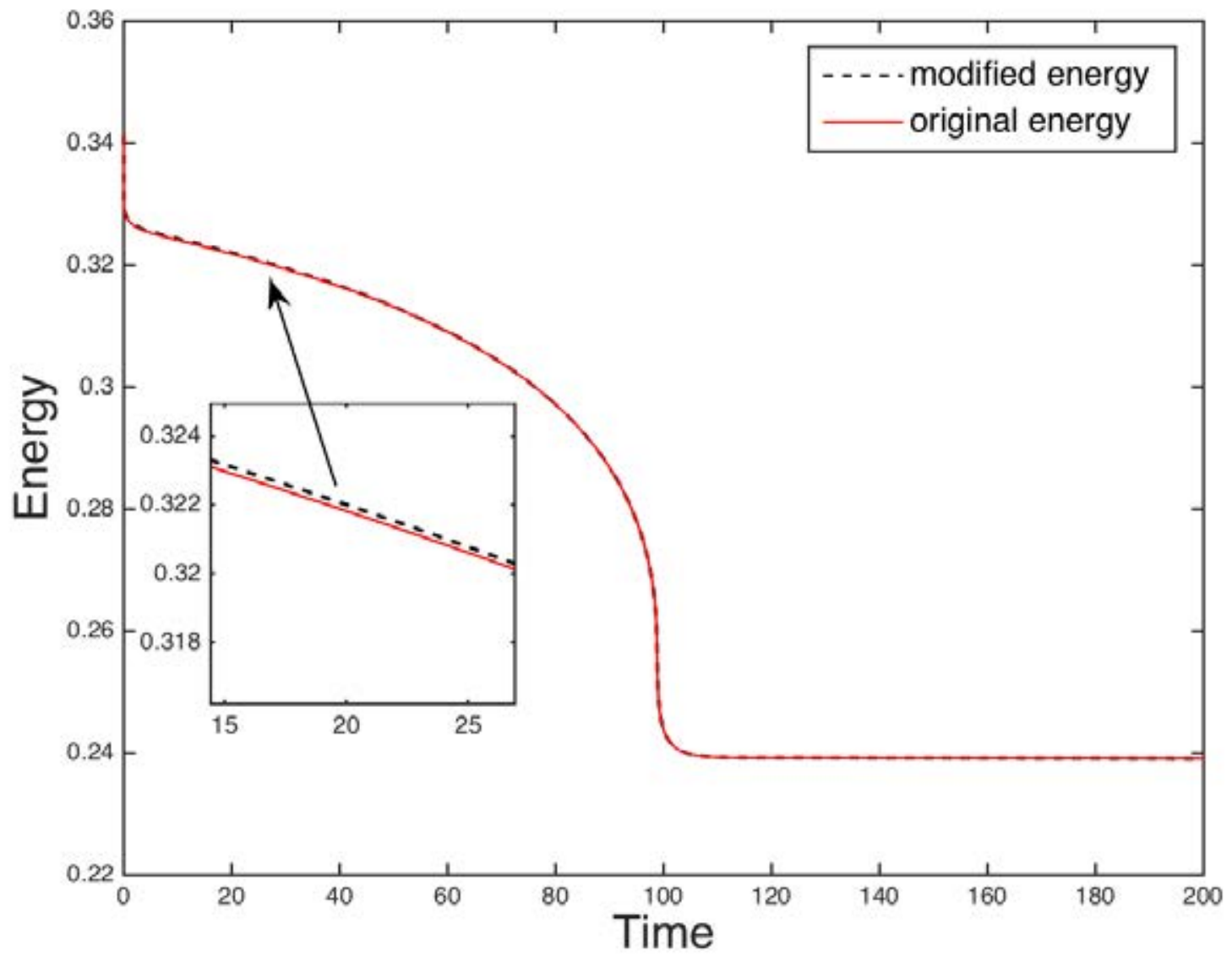

FIGURE. 4.5. Time evolution of the two free energy functionals for the Example 4.2, the modified energy (3.47) and the original energy (3.2), by using the initial condition (4.3). The small inset figure is a close-up view.

4.3. **Spinodal decomposition in 2D and 3D.** In this example, we study the phase separation behaviors in the 2D and 3D spaces that are called spinodal decomposition. The process of the phase separation can be studied by considering a homogeneous binary mixture, which is quenched into the unstable part of its miscibility gap. In this case, the spinodal decomposition takes place, which manifests in the spontaneous growth of the concentration fluctuations that leads the system from the homogeneous to the two-phase state. Shortly after the phase separation starts, the domains of the binary components are formed and the interface between the two phases can be specified.

The initial conditions are taken as randomly perturbed fields. For 2D, it reads as

$$\phi_0(x,y) = 0.4 + 0.001\text{rand}(x,y), \quad \rho_0(x,y) = 0.3; \tag{4.4}$$

and for 3D, it reads as

$$\phi_0(x,y,z) = 0.4 + 0.001\text{rand}(x,y,z), \quad \rho_0(x,y) = 0.3, \tag{4.5}$$

where the rand$(\cdot)$ is the random number in $[-1,1]$ and has zero mean.

We choose the time step size $\delta t = 1\text{e}{-3}$ to perform both of the 2D and 3D simulations. In Fig. 4.6, we show the snapshots of $\phi$ and $\rho$ for the 2D simulations, where we observe the coarsening dynamics that the fluid component with the less volume accumulates to small satellite drops everywhere. When the time evolves, the small drops will collide, merge and form drops with bigger sizes. The final equilibrium solution is obtained around after $t = 1700$, where all small bubbles accumulate into a big bubble.

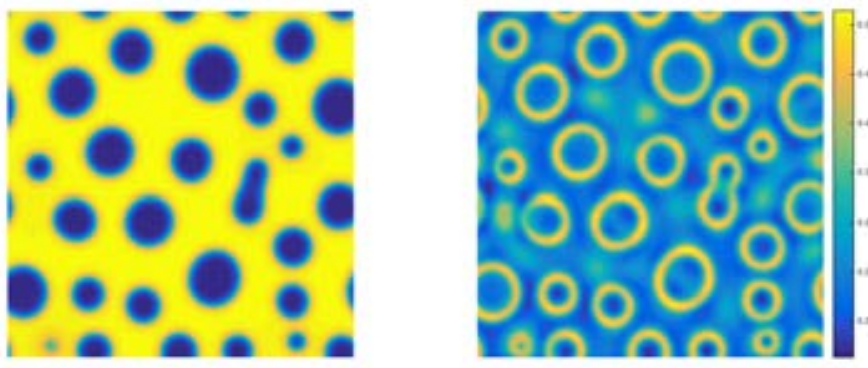

(a) $t = 10$.

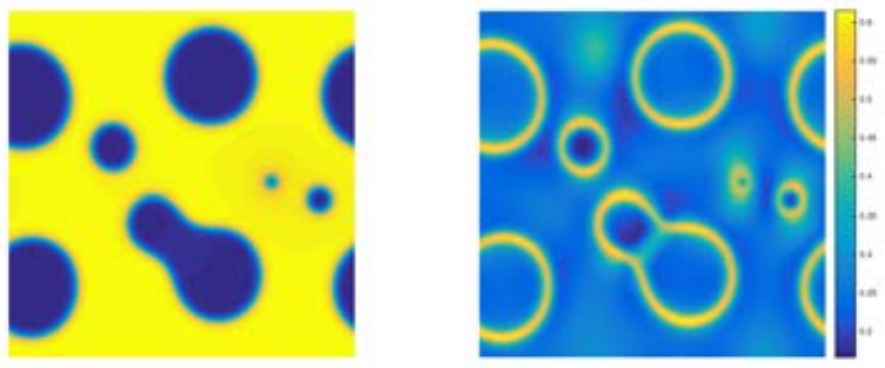

(b) $t = 80$.

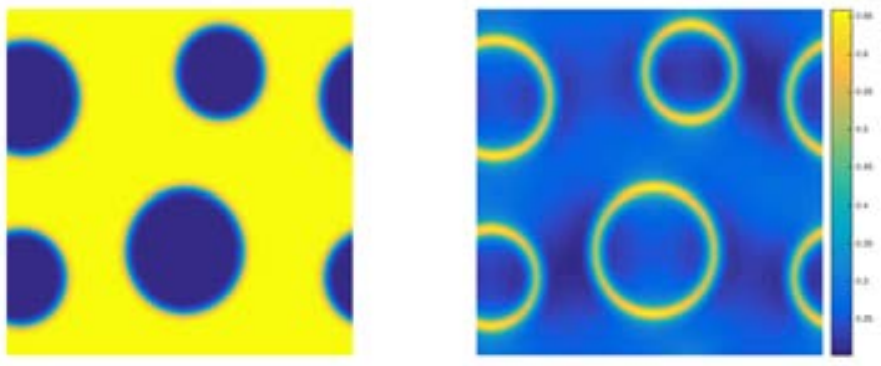

(c) $t = 200$.

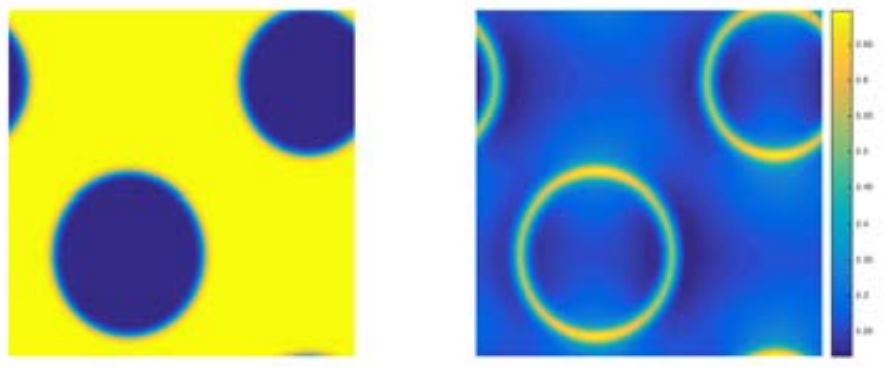

(d) $t = 500$.

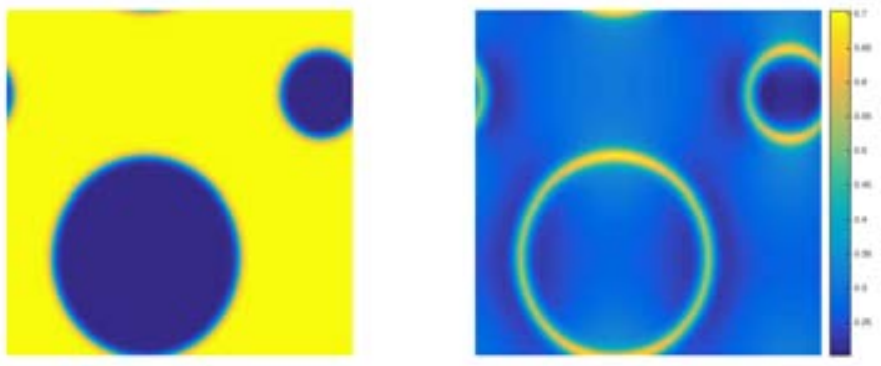

(e) $t = 1500$.

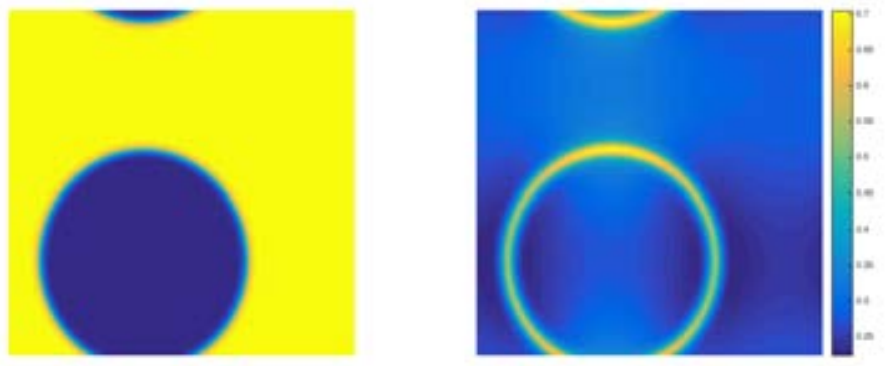

(f) $t = 2000$.

FIGURE. 4.6. 2D spinodal decomposition for random initial data of (4.4). Snapshots of phase variables $\phi$ and $\rho$ are taken at $t = 10$, 80, 200, 500, 1500 and 2000. For each panel, the left subfigure is the profile of $\phi$, and the right subfigure is the profile of $\rho$.

In Fig. 4.7, we show the snapshots of $\phi$ and $\rho$ for the 3D simulations. Similar to the 2D case, the two fluids initially are well mixed, and they sooner start to decompose and accumulate. The final steady state (shown in Fig. 4.7(d)) is consistent with the 2D results, namely, all small bubbles accumulate to form the final big bubble. In order to obtain more accurate view, since the computed domain is periodic, in Fig. 4.8, we plot the isosurface for 4 periods, i.e., $[0, 8\pi]^3$. In Fig. 4.9, we present the evolution of total free energy (3.47) for both of these two simulations, that clearly shows the energies always monotonically decay with respect to the time.

### 4.4. **Surfactant absorption.**

4.4.1. *Surfactant uniformly distributed initially.* We assume the fluid interface and the surfactant are uniformly distributed over the domain initially and the specific profiles (shown in Fig. 4.10 (a)) are

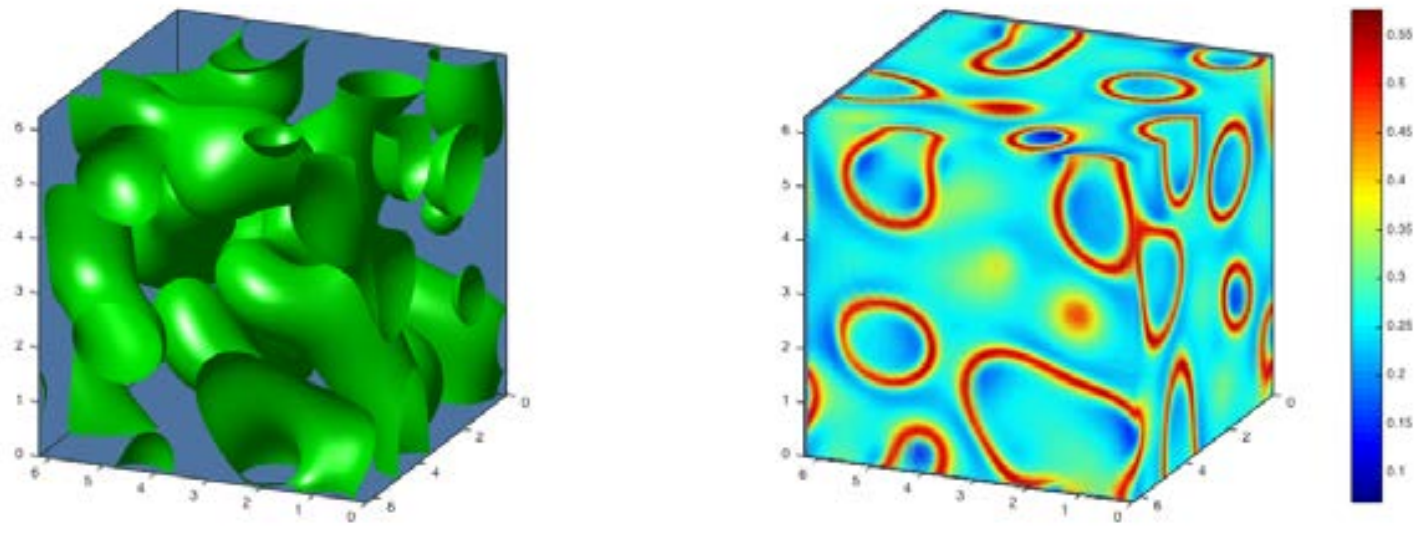

(a) $t = 5$.

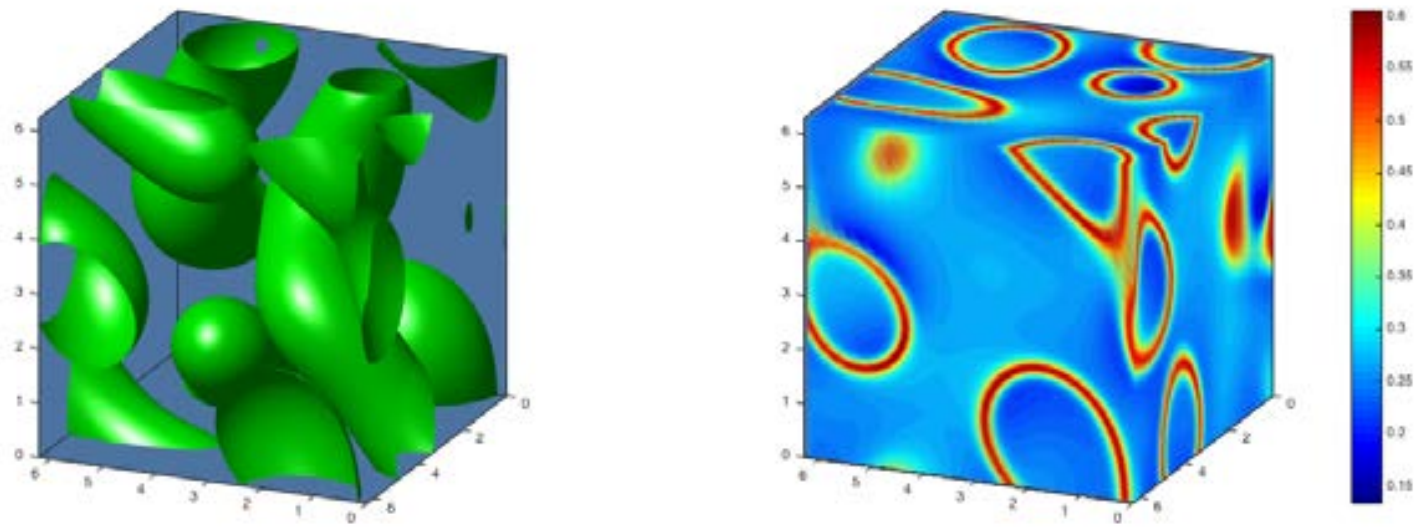

(b) $t = 15$.

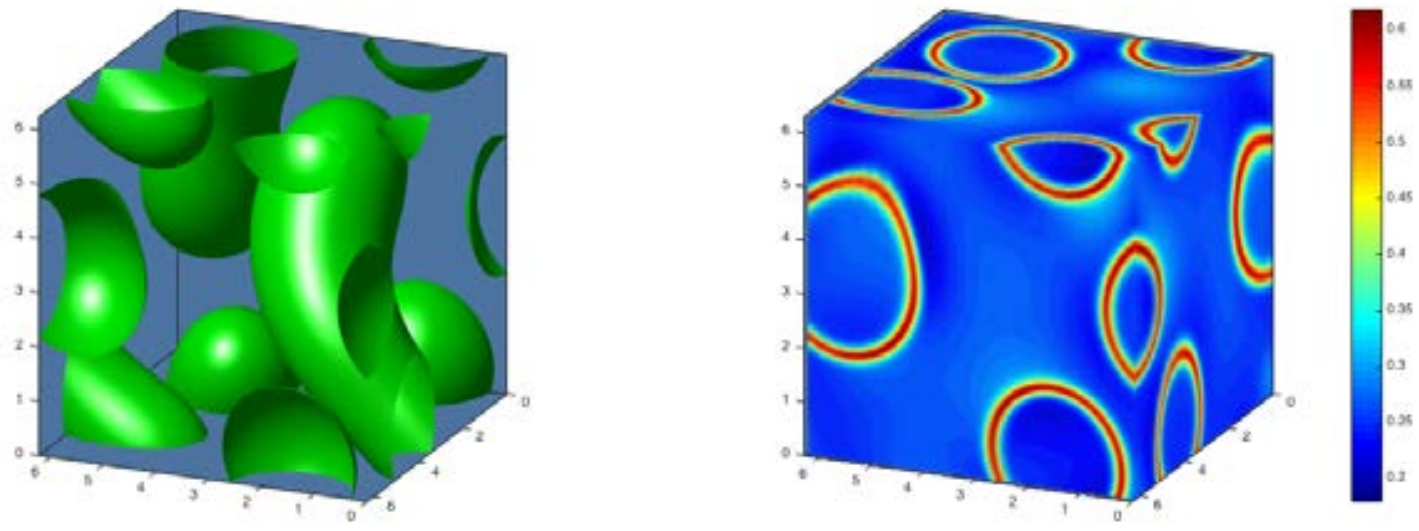

(c) $t = 25$.

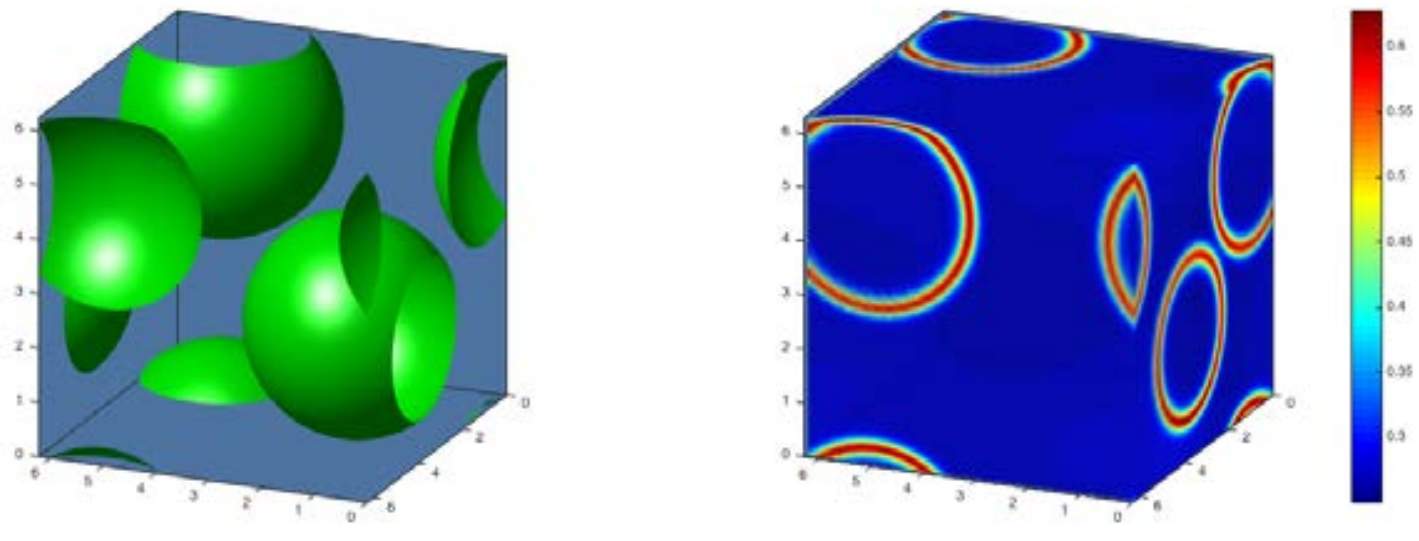

(d) $t = 95$.

FIGURE. 4.7. 3D spinodal decomposition for random initial data of (4.5). Snapshots of phase variables $\phi$ and $\rho$ are taken at $t = 5$, 15, 25, and 95. For each panel, the left subfigure is the profile of $\phi$ (isosurfaces of $\{\phi = 0\}$), and the right subfigure is the profile of $\rho$ (that is visualized by 128 cross-sections along the z-axis).

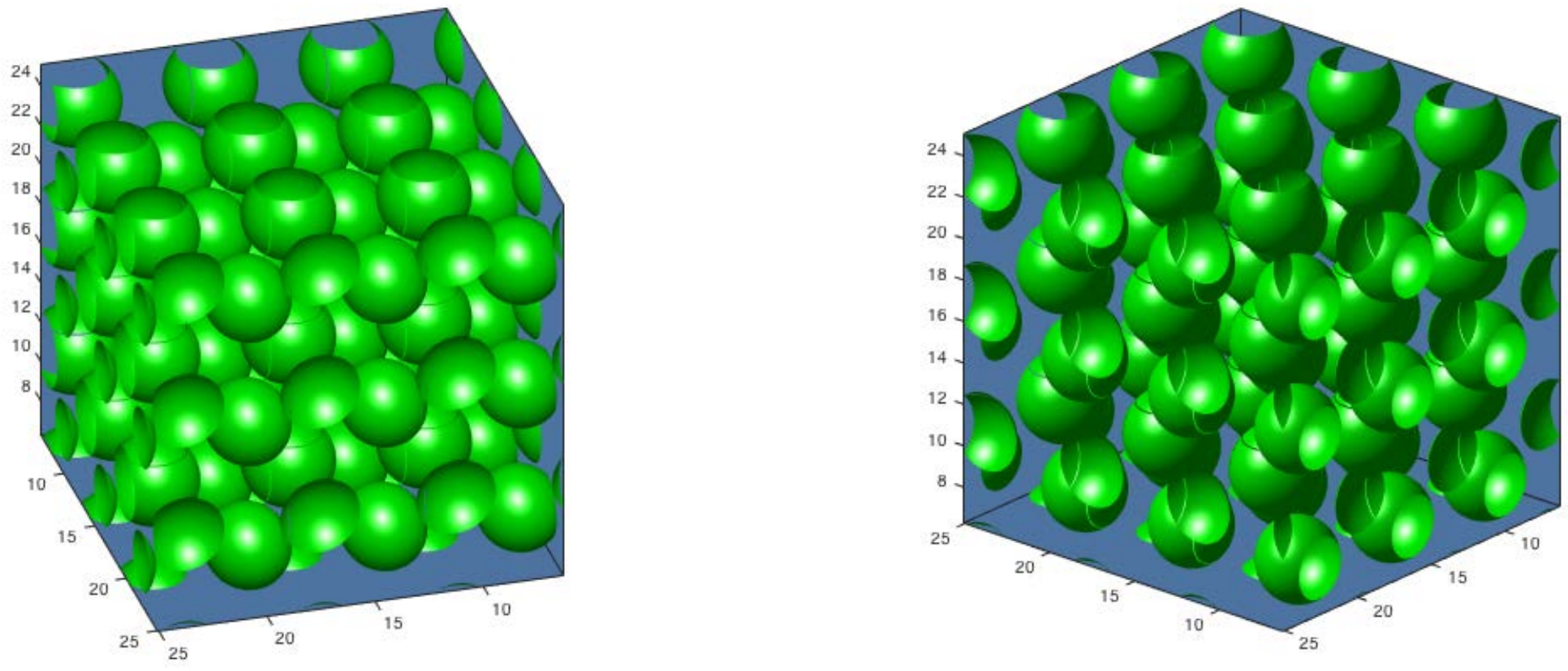

FIGURE. 4.8. The isosurfaces $\{\phi = 0\}$ of the steady state solution at $t = 95$ of the 3D spinodal decomposition example for 4 periods, i.e., $[0, 8\pi]^3$. The two subfigures are from different angles of view.

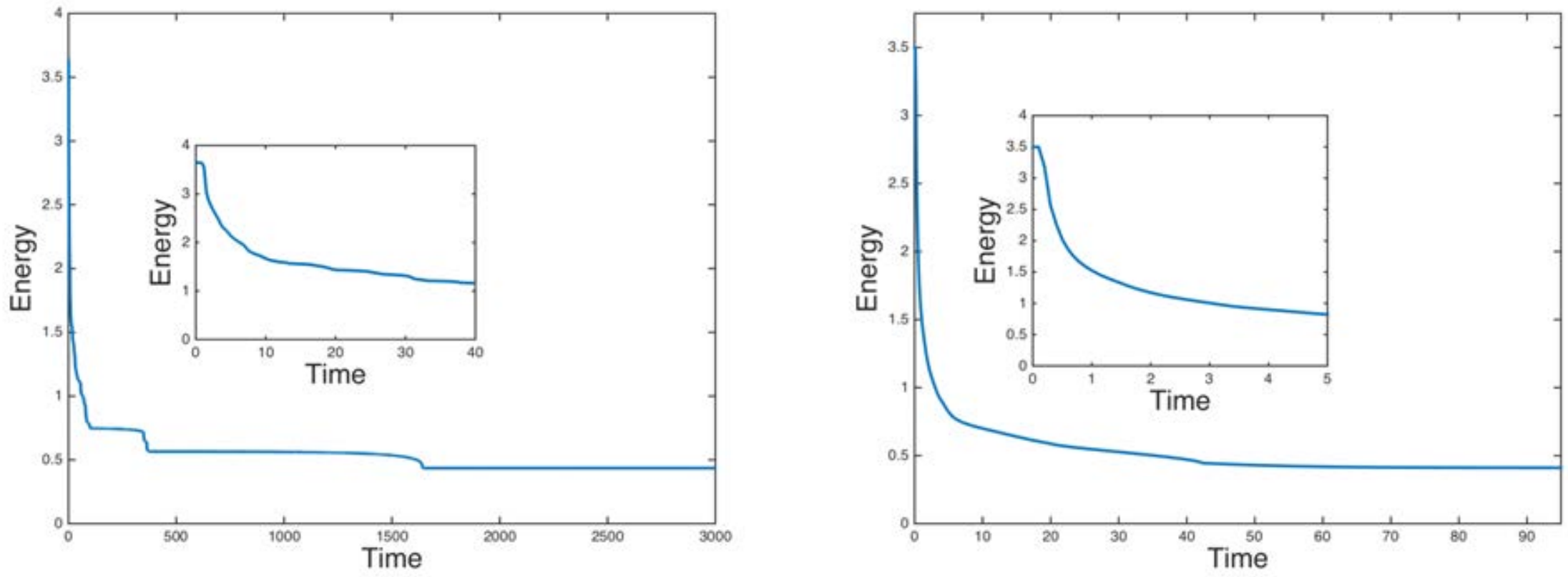

FIGURE. 4.9. (From left to right:) Time evolution of the free energy functional (3.47) of the 2D and 3D spinodal decompositions, respectively. The small inset figure is a close-up showing where the energy decreases fast.

chosen

$$\phi_0(x, y) = 0.1 + 0.01\cos(6x)\cos(6y), \tag{4.6}$$

$$\rho_0(x, y) = 0.2 + 0.01\cos(6x)\cos(6y). \tag{4.7}$$

We take the time step size $\delta t = 1\text{e}{-3}$. Fig. 4.10 shows the snapshots of coarsening dynamics at $t = 0, 3, 30, 100, 300, 1500$ and the final steady shape forms a big drop. Driven by the coupling entropy energy term, the surfactant is absorbed into the binary fluid interfaces so that the higher concentration appears near the interfaces than other regions. In Fig. 4.11, we show the profiles of the velocity $\boldsymbol{u}$ and $p$ at $t = 30$, 300 and 1500.

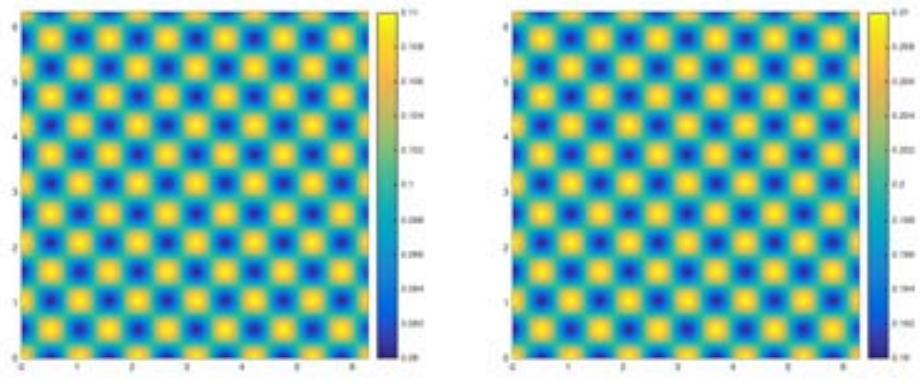

(a) $t = 0$.

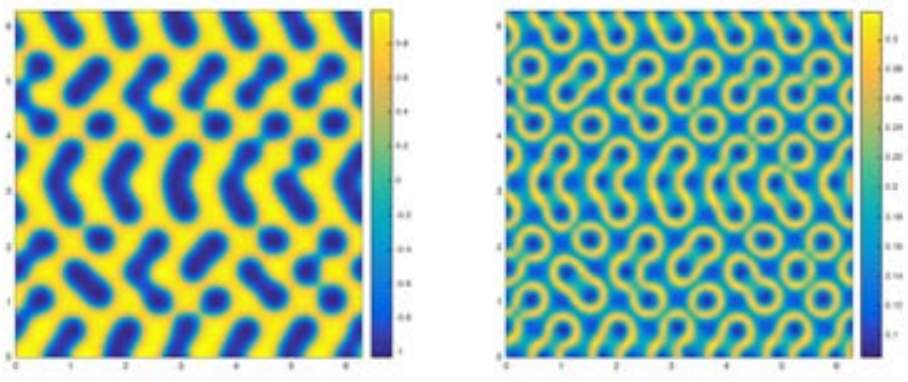

(b) $t = 3$.

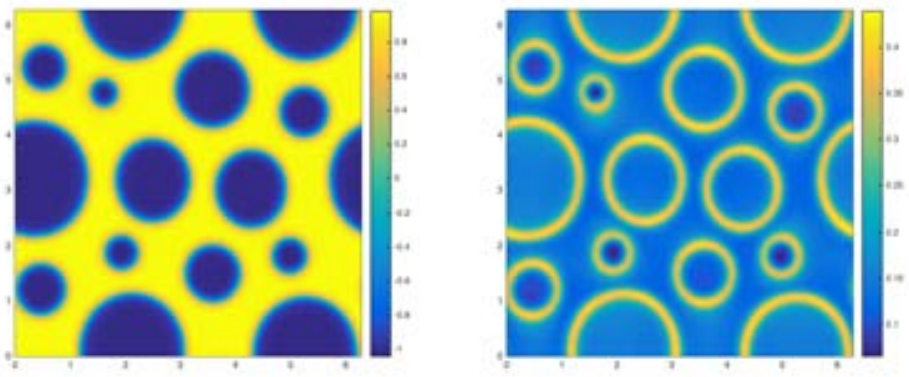

(c) $t = 30$.

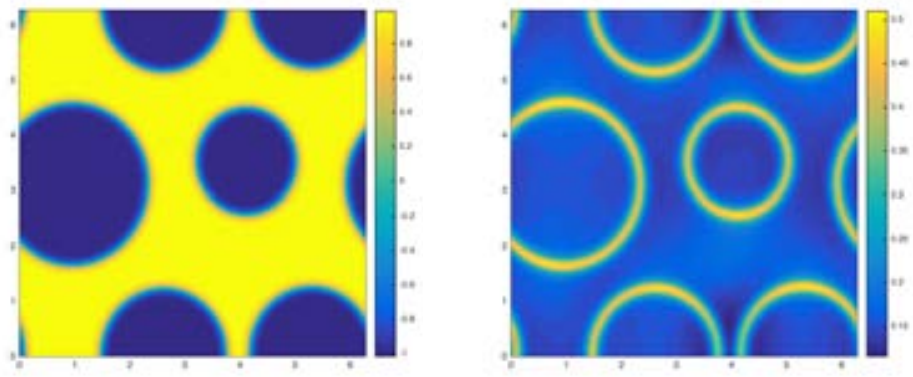

(d) $t = 100$.

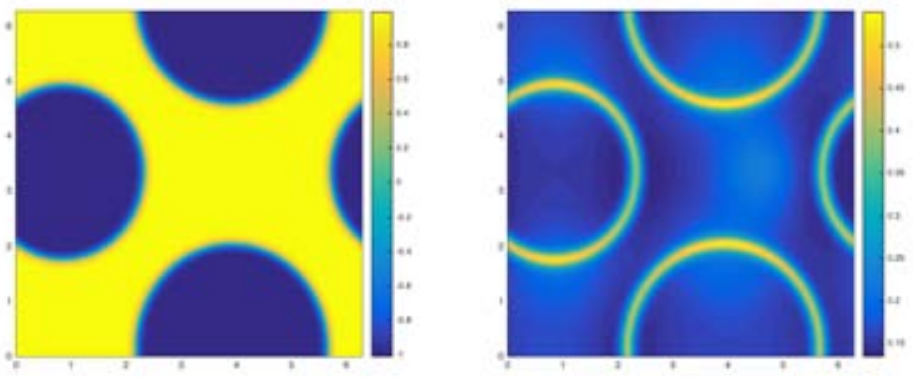

(e) $t = 300$.

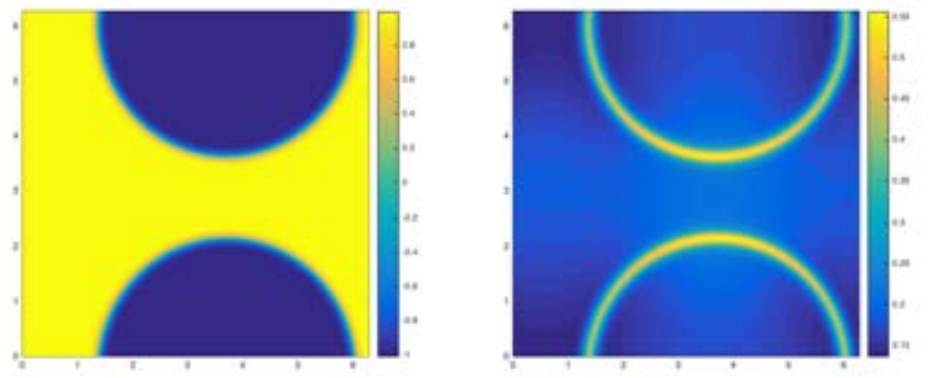

(f) $t = 1500$.

FIGURE. 4.10. Snapshots of the phase variables $\phi$ and $\rho$ are taken at $t = 0$, 3, 30, 100, 300, and 1500 for Example 4.4.2, where the surfactants are distributed uniformly initially. For each panel, the left subfigure is the profile of $\phi$ and the right subfigure is the profile of $\rho$.

4.4.2. *Surfactant locally distributed initially.* We here assume the fluid interface and the surfactant field are mismatched over the domain initially. When $t = 0$, the phase field variable $\phi$ is same as the previous example, but the surfactant concentration variable $\phi$ is locally accumulated around the center of the domain. The specific initial profiles (shown as Fig. 4.12 (a)) are chosen as

$$\phi_0(x, y) = 0.1 + 0.01\cos(6x)\cos(6y), \tag{4.8}$$

$$\rho_0(x, y) = 0.8\exp\Big(-\frac{(x-\pi)^2 + (y-\pi)^2}{1.25^2}\Big). \tag{4.9}$$

We take the time step size $\delta t = 1\mathrm{e}{-3}$ due to the better accuracy. Fig. 4.12 shows the snapshots of the dynamics behaviors of $\phi$ and $\rho$ at various time. Since the surfactant is initially concentrated at the center, it takes longer times for the surfactant to diffuse away from this center region. Consequently, during the early stage of the evolution, the higher concentration of surfactant only appears around the center area of the domain. We observe that, the surfactant totally completely diffuses and is absorbed into the binary fluid interfaces after around $t = 300$.

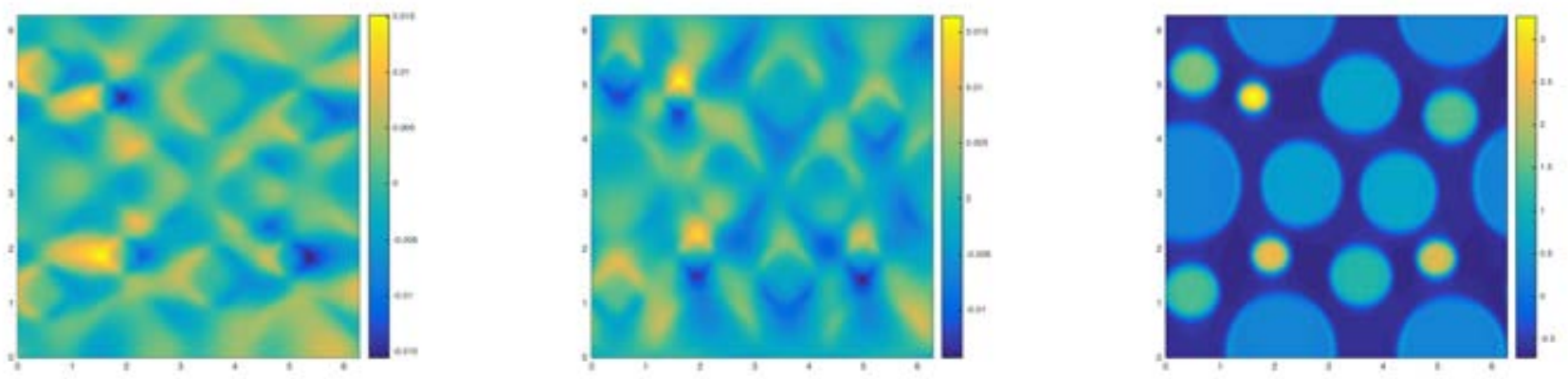

(a) The first component of the velocity field $\boldsymbol{u} = (u, v)$ at $t = 30, 300, 1500$.

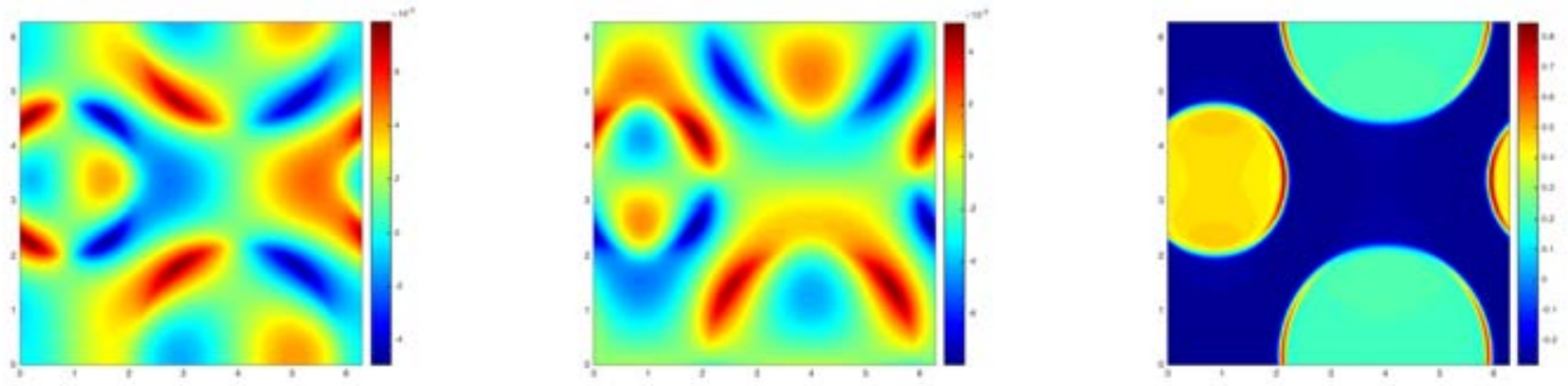

(b) The second component of the velocity field $\boldsymbol{u} = (u, v)$ at $t = 30, 300, 1500$.

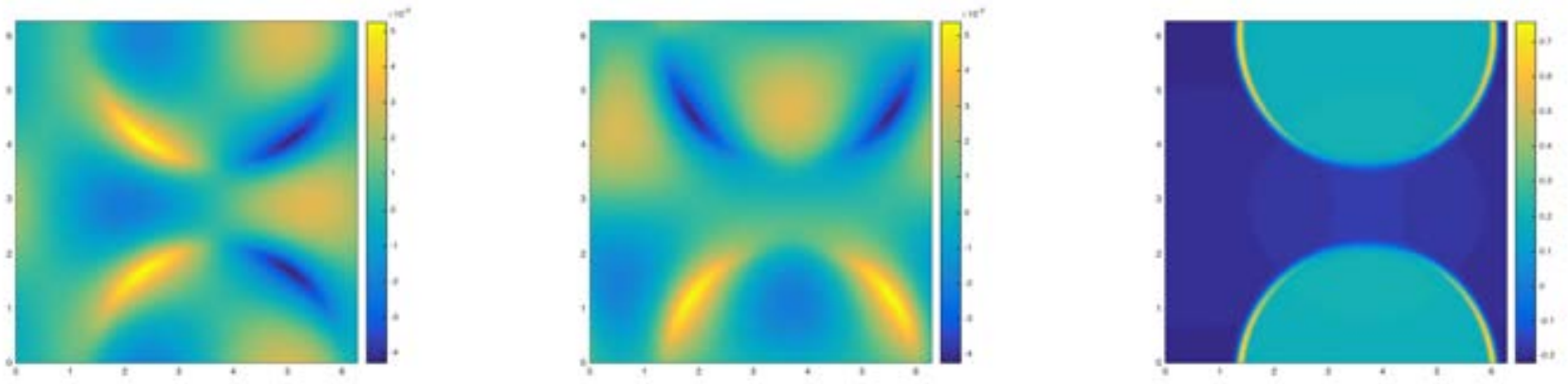

(c) The pressure $p$ at $t = 30, 300, 1500$.

Figure. 4.11. Profiles of the velocity field $\boldsymbol{u} = (u, v)$ and the pressure $p$. Snapshots are taken at $t = 30$, 300 and 1500 for the Example 4.4.1.

We finally plot the evolution of energy curves in Fig. 4.13 for the Example 4.4.1 and the Example 4.4.2, respectively. Overall, the numerical solutions of both of the two examples present similar features to those obtained in [35, 11].

## 5. Concluding remarks

In this paper, we have developed two efficient schemes for solving the binary fluid-surfactant phase field model coupled with the fluid flow. We combine the projection method, the IEQ approach, and implicit-explicit treatments for the stress and convective terms. Our schemes are purely linear, second order accurate in time. Furthermore, the induced linear systems are well-posed and the schemes are provably unconditionally energy stable. The developed schemes can serve as a building block to design accurate and energy stable linear schemes for fluid flow coupled gradient flow problems with single or multiple variables. Although we consider only time discrete schemes in this paper, the results here can be carried over to any consistent finite-dimensional Galerkin type approximations since the proofs are

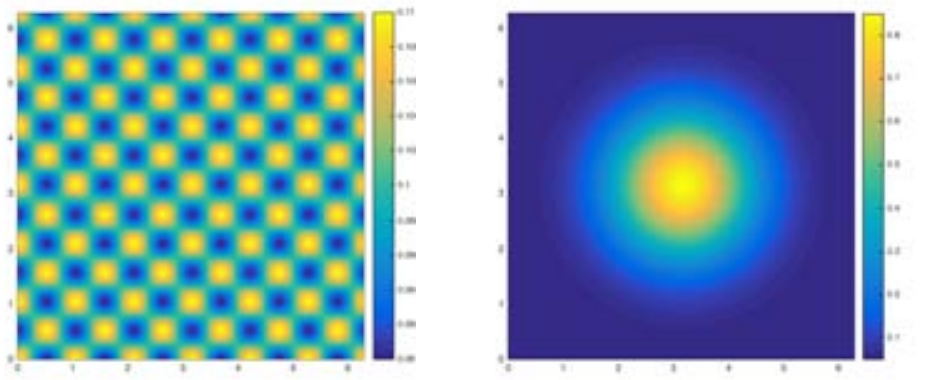

(a) $t = 0$.

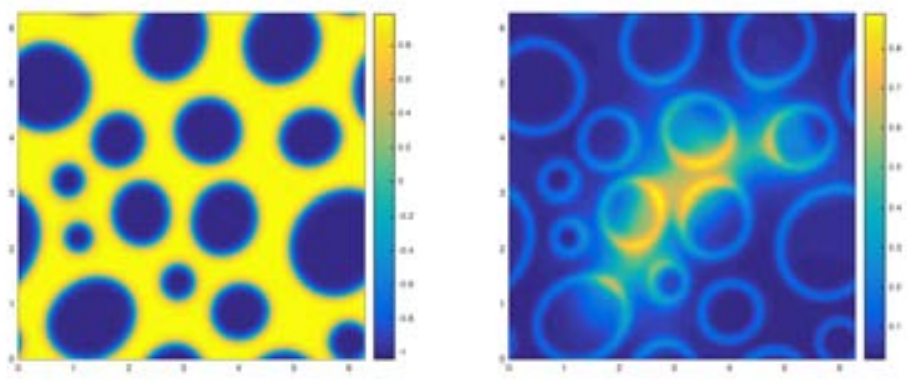

(b) $t = 10$.

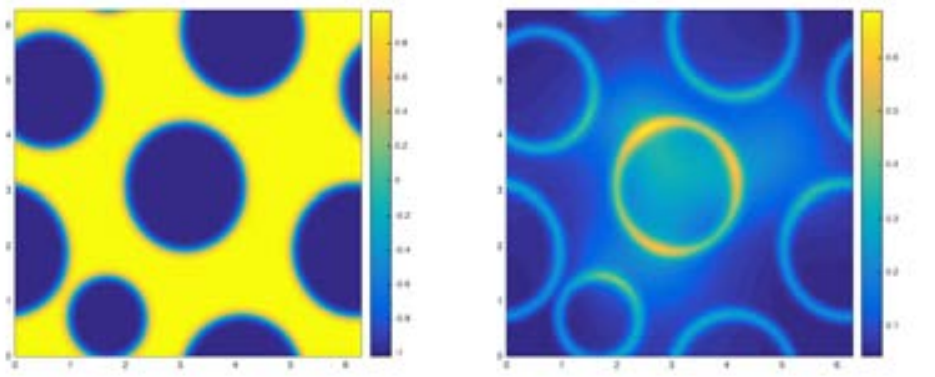

(c) $t = 100$.

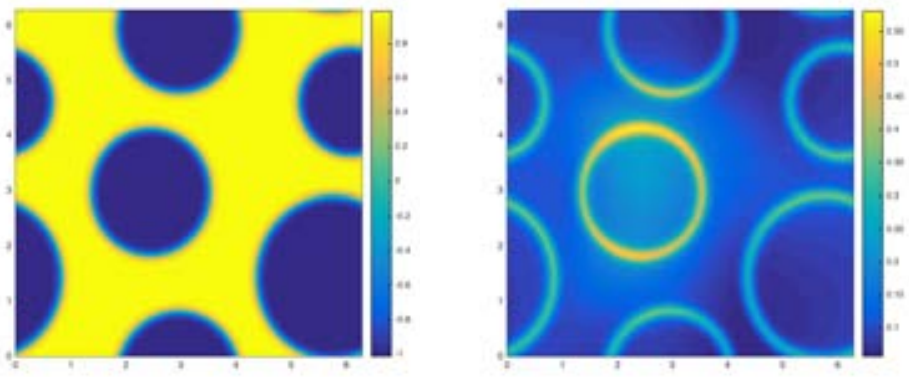

(d) $t = 200$.

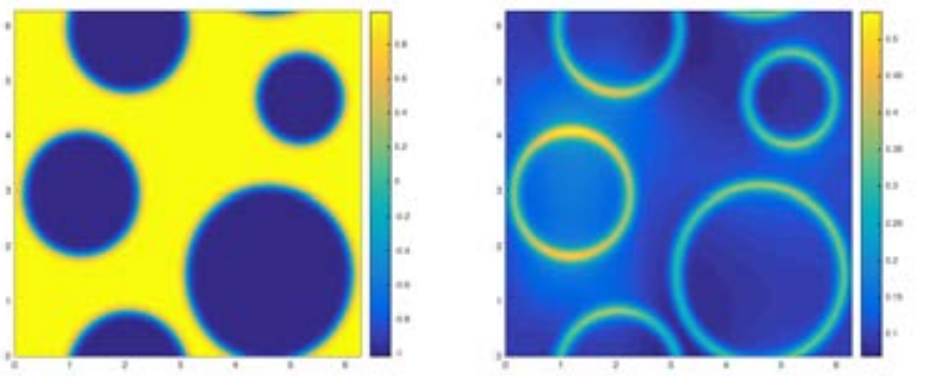

(e) $t = 300$.

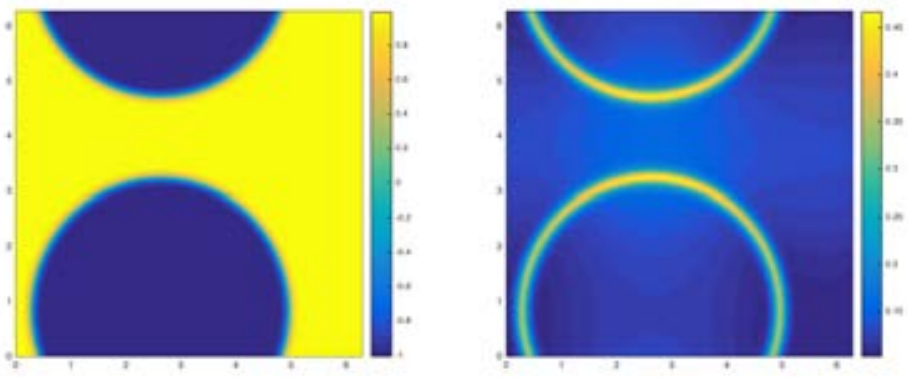

(f) $t = 1000$.

FIGURE. 4.12. Snapshots of the phase variables $\phi$ and $\rho$ are taken at $t = 0$, 10, 100, 200, 300, and 1000 for Example 4.4.2, where the surfactants are initially distributed in a center circle. For each panel, the left subfigure is the profile of $\phi$ and the right subfigure is the profile of $\rho$.

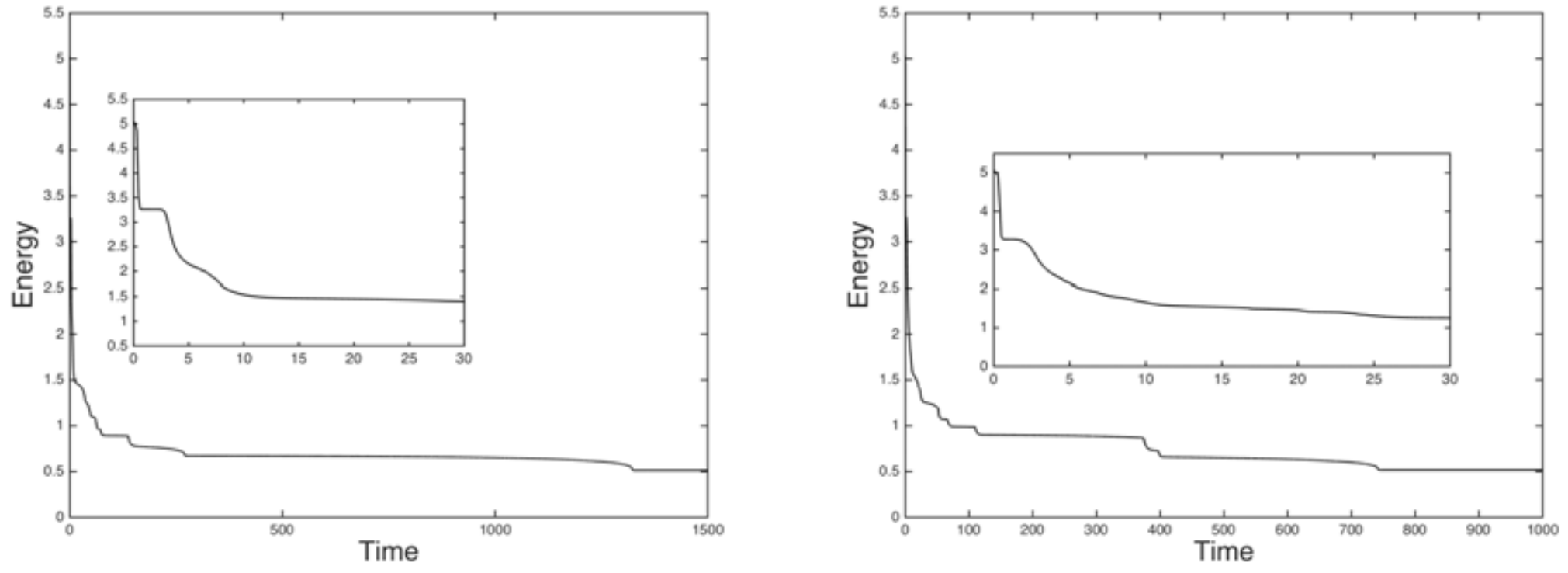

FIGURE. 4.13. Time evolution of the free energy functional (3.6) for the Examples 4.4.1 and 4.4.2, respectively. The small inset figure is a close-up view showing where the energy decreases fast.

all based on a variational formulation with all test functions in the same space as the space of the trial functions.